\documentclass[reqno]{amsart}
\usepackage{amsfonts}
\usepackage{graphicx}
\usepackage{amsmath,amssymb,amsfonts,amsthm,mathrsfs}
\usepackage[english]{babel}
\usepackage{hyperref}
\usepackage{marginnote}
\usepackage{color}
\usepackage{enumerate}

\newcommand{\dis}{\displaystyle}

\newcommand{\R}{\mathbb{R}}

\theoremstyle{plain}

\newtheorem{theorem}{Theorem}[section]
\newtheorem{lemma}[theorem]{Lemma}

\newtheorem{definition}[theorem]{Definition}
\newtheorem{remark}[theorem]{Remark}
\theoremstyle{remark}
\numberwithin{equation}{section}

\def\v{\varepsilon}

\def\t{\theta}

\def\a{\alpha}

\def\g{\gamma}
\def\d{\delta}

\def\s{\sigma}

\def\f{\frac}

\usepackage{tikz}

\newcommand{\dd}{{\rm d}}

\newcommand{\FL}{\mathbf{L}}

\newcommand{\Fi}{\mathbf{1}}
\newcommand{\la}{\lambda}

\newcommand{\vep}{{\varepsilon}}
\newcommand{\fu}{\mathfrak{u}}

\begin{document}

\title[Boundary layer for specular boundary condition]{Boundary layer solution of the Boltzmann equation for specular boundary condition}

\author[F.M. Huang]{Feimin Huang}
\address[F. M. Huang]{Academy of Mathematics and Systems Science, Chinese Academy of Sciences, Beijing 100190, China; School of Mathematical Sciences, University of Chinese Academy of Sciences, Beijing 100049, China
}
\email{fhuang@amt.ac.cn}

\author[Z.H. Jiang]{Zaihong Jiang}
\address[Z.H. Jiang]{ Department of Mathematics, Zhejiang Normal University, Jinhua 321004, China
}
\email{jzhong@zjnu.cn}

\author[Y. Wang]{Yong Wang}
\address[Y. Wang]{Academy of Mathematics and Systems Science, Chinese Academy of Sciences, Beijing 100190, China; School of Mathematical Sciences, University of Chinese Academy of Sciences, Beijing 100049, China
}
\email{yongwang@amss.ac.cn}

\begin{abstract}
In the paper, we establish the existence of steady boundary layer solution of Boltzmann equation with specular boundary condition in $L^2_{x,v}\cap L^\infty_{x,v}$ in half-space. The uniqueness, continuity and exponential decay of the solution are obtained, and such estimates are important to prove the Hilbert expansion of Boltzmann equation for half-space problem with specular boundary condition.
\end{abstract}

\date{\today}
\maketitle

\tableofcontents

\thispagestyle{empty}


\section{Introduction}

In the present paper, we consider the steady Boltzmann equation
\begin{equation}\label{1.1}
v_3\cdot \partial_x F=Q(F,F)+\mathcal{S},\quad \   (x,v)\in \mathbb{R}_+\times\mathbb{R}^3,
\end{equation}
with $\R_+=(0,+\infty)$.
The Boltzmann collision term $Q(F_1,F_2)$ on the right is defined in terms of the following bilinear form
\begin{align}\label{1.2}
Q(F_1,F_2)&\equiv\int_{\mathbb{R}^3}\int_{\mathbb{S}^2} B(v-u,\t)F_1(u')F_2(v')\,{d\omega du}
-\int_{\mathbb{R}^3}\int_{\mathbb{S}^2} B(v-u,\t)F_1(u)F_2(v)\,{d\omega du}\nonumber\\
&:=Q_+(F_1,F_2)-Q_-(F_1,F_2),
\end{align}
where the relationship between the post-collision velocity $(v',u')$ of two particles with the pre-collision velocity $(v,u)$ is given by
\begin{equation*}
u'=u+[(v-u)\cdot\omega]\omega,\quad v'=v-[(v-u)\cdot\omega]\omega,
\end{equation*}
for $\omega\in \mathbb{S}^2$, which can be determined by conservation laws of momentum and energy
\begin{equation*}
u'+v'=u+v,\quad |u'|^2+|v'|^2=|u|^2+|v|^2.
\end{equation*}
The Boltzmann collision kernel $B=B(v-u,\theta)$ in \eqref{1.2} depends only on $|v-u|$ and $\theta$ with $\cos\theta=(v-u)\cdot \omega/|v-u|$.   Throughout this paper,  we  consider the hard sphere model, i.e.,
\begin{equation*}
B(v-u,\t)=|(v-u)\cdot \omega|
\end{equation*}
We supplement the Boltzmann equation \eqref{1.1} with the perturbed specular reflection boundary condition
\begin{align}\label{B.C}
F(0,v)|_{v_3>0}=F(0,Rv)+F_b(Rv)
\end{align}
where $Rv=(v_1,v_2,-v_3)$, and $F_b(v)$ is a given function. We impose the condition at infinity
\begin{equation}\label{Farfield}
\lim_{x\rightarrow\infty}F(x,v)=\mu\equiv \f{1}{(2\pi)^{\frac32}}\exp\{-\f{|v-\mathfrak{u}|^2}{2}\},
\end{equation}
where $\mathfrak{u}=(\fu_1,\fu_2,\fu_3)\in\mathbb{R}^3$ is a given back ground macroscopic velocity which is independent of $x$. Throughout the paper, we always assume  $\fu_3=0$ which insure that $\mu(v)=\mu(R v)$, i.e., $\mu$ satisfies the specular boundary condition.

\smallskip

We look for solutions in the form
\begin{equation*}
\mathfrak{f}(x,v)=\frac{F(x,v)-\mu}{\sqrt{\mu}},
\end{equation*}
then \eqref{1.1},\eqref{B.C} and \eqref{Farfield} are  rewritten as
\begin{align}\label{1.7}
\begin{cases}
v_3 \partial_x \mathfrak{f}+\FL \mathfrak{f}=\Gamma(\mathfrak{f},\mathfrak{f})+S,\\
\mathfrak{f}(0,v)|_{v_3>0}=\mathfrak{f}(0,Rv)+f_b(Rv),\\
\lim_{x\rightarrow\infty}\mathfrak{f}(x,v)=0.
\end{cases}
\end{align}
where we have denoted
\begin{align*}
\begin{split}
&\FL \mathfrak{f}=-\frac1{\sqrt\mu} \{Q(\mu,\sqrt{\mu} \mathfrak{f})+Q(\sqrt{\mu} \mathfrak{f},\mu)\},\\
&\Gamma (\mathfrak{f},\mathfrak{f})=\frac1{\sqrt\mu} Q(\sqrt{\mu} \mathfrak{f}, \sqrt{\mu} \mathfrak{f}),\\
&f_b(v)=\frac{F_b(v)}{\sqrt{\mu}},\quad S= \frac{\mathcal{S}}{\sqrt{\mu}}.
\end{split}
\end{align*}
It is noted that the function $f_b(v)$ is defined only for $v_3< 0$, and we  assume that it is extended to be $0$ for $v_3\geq0$ throughout the paper.

\smallskip

There have been many studies on the half-space problem of the steady Boltzmann equation in the literature.  The existence, uniqueness and properties of asymptotic behavior were proved in \cite{BCN} for the linearized Boltzmann equation of a hard sphere gas for the  Dirichlet type boundary condition, see \cite{CGS} for a classification of well-posed kinetic boundary layer problem. Later, the existence of nonlinear Knudsen boundary layers with small magnitudes for the hard sphere model was established in \cite{UYU}, it was shown that the existence of a solution depends on the Mach number of the far field Maxwellian, and an implicit solvability conditions yielding the co-dimensions of the boundary data, see \cite{UYU1,DWY} for the time-asymptotic stability of such boundary layer solution; we also refer \cite{Wu} for the construction of a modified boundary layer solution in $L^\infty_{x,v}$ space which is used to prove the Hilbert expansion in a disk;  recently,  for the purpose of studying the transition from evaporation to condensation,  Bernhoff-Golse \cite{BG1} offered the existence and uniqueness of a uniformly decaying boundary layer type solution in the situation that gas is in contact with its condensed phase.
For the diffuse boundary condition, the existence of steady Boltzmann solution is proved in \cite{EGKM,DHWZ} in boundary domain, and the time-asymptotic stability of such steady solutions is also obtained.
For  the specular reflection condition and the solution tends to a global Maxwellian in the far field, Golse-Perthame-Sulem \cite{GPS-1988} proved the existence, uniqueness and asymptotic behavior  in the functional space \eqref{1.08}.
To prove the Hilbert expansion of Boltzmann equation for half-space problem with specular boundary condition, the continuity and uniform estimate in $L^\infty_{x,v}$ are needed,  so in the present paper, we will focus on the existence steady solution of \eqref{1.1}  in the functional space $L^{2}_{x,v}\cap L^\infty_{x,v}$.

\smallskip

Now we list some notations that will be used in this paper.
Throughout this paper, $C$ denotes a generic positive constant which may vary from line to line.
And $C_a$, $C_b, \cdots$ denote the generic positive constants depending on $a, b, \cdots$, respectively, which also may vary from line to line. $\| \cdot \|_{L^p}$  denotes the standard $L^p(\Omega \times \mathbb{R}_v^3)$-norm or $L^p(\mathbb{R}_v^3)$-norm, $\|\cdot \|_{\nu}=\|\sqrt{\nu}\cdot \|_{L^2}$. When the norms need to be distinguished from each other, we write $\|\cdot\|_{L^p_v}$,
$\|\cdot\|_{L^p_x}$ and $\|\cdot\|_{L^p_{x,v}}$
respectively.  Moreover, for the phase boundary integration, we define $|\cdot|_{L^\infty(\gamma)}$ denotes the $L^\infty(\gamma)$-norm,
$|\cdot|_{L^2(\gamma)}$ denotes the $L^2(\gamma, |v_3|dv)$-norm,
where $\gamma=\partial \Omega \times  \mathbb{R}_v^3 $ with $\Omega=\R_+\ \mbox{or} \  (0,d)$ for $d>0$.

\

We define the weight function
\begin{align}\label{wt}
w(v)=(1+|v|^2)^{\frac{\beta}{2}}e^{\varsigma |v-\mathfrak{u}|^2}.
\end{align}
where $\beta, \varsigma $ are two positive constants.
\begin{theorem}	\label{thm1.1}
Let $\beta \geq 3$ and $0\leq \varsigma<\f14$. We assume $S\in \mathbb{N}^{\perp}$ and $f_b(v)$ satisfies
\begin{align}\label{1.12}
\int_{\R^3}f_b(v) v_3  \sqrt{\mu}dv
&=\int_{\R^3}(v_1-\fu_1)v_3 f_b(v) \sqrt{\mu}dv
=\int_{\R^3}(v_2-\fu_2)v_3 f_b(v) \sqrt{\mu}dv\nonumber\\
&=\int_{\R^3}v_3 |v-\fu|^2f_b(v) \sqrt{\mu}dv=0,
\end{align}
There exists a small $\delta_0>0$ such that if $| wf_b |_{L_v^\infty}+\|\nu^{-1} w  e^{\sigma_0x}S\|_{L^\infty_{x,v}}\leq \d_0$,  then the boundary value problem \eqref{1.7} has a unique solution satisfying
\begin{align}\label{ES}
\|e^{\sigma x}w f\|_{L^\infty_{x,v}}
&\leq \frac{C}{\sigma_0-\sigma} (| w  f_b |_{L^\infty_v}+\|\nu^{-1}w  e^{\sigma_0x}S\|_{L^\infty_{x,v}}),
\end{align}
where $\sigma $ is a constant  such that $0<\sigma< \sigma_0$, and $C>0$ is a constant independent of $\sigma$,  and $\delta_0$ depends on $\frac{1}{\sigma_0-\sigma}$ with $\delta_0\rightarrow 0+$ as $\sigma\rightarrow \sigma_0$.
Moreover, if
$S$ is continuous in $\R_+\times \mathbb{R}^3$ and $f_b$ is continuous in $\{v  \in\R^3\}$, then $f(x,v)$ is continuous away from the  grazing set $\big\{(0,v)\ : \ v\in\R^3,\  v_3=0\big\}$.
\end{theorem}

\begin{remark}
Golse-Perthame-Sulem \cite{GPS-1988} proved an existence result for \eqref{1.7} in the following functional space
\begin{equation}\label{1.08}
\int_{\R_+}\int_{\R^3} e^{\s x} (1+|v|) |f(x,v)|^2dvdx+\int_{\R^3} \sup_{x\in\R_+}\{|e^{\s x} f(x,v)|^2\} dv<\infty.
\end{equation}
In the present paper, we are interested in the existence result of Boltzmann equation \eqref{1.1}  in the functional space  $L^{2}_{x,v}\cap L^\infty_{x,v}$ with uniqueness, continuity;  in particular, we shall use the linear existence result in Theorem \ref{thm3.1} below to prove the Hilbert expansion of Boltzmann equation for half-space problem in \cite{GHW}.
\end{remark}

 This paper is organized as follows: In section \ref{P}, we present some useful results which will be used later. In section \ref{L}, we study the existence theory of the linearized problem with a source term. The proof of our main theorem \ref{thm1.1} is given in section \ref{M}.

\

\section{Preliminary}\label{P}

As in \cite{Gl}, $\FL$ can be decomposed as $\FL=\nu(v)-K$ where
\begin{align}\notag
(K_1f)(v)&=\int_{\mathbb{R}^3}\int_{\mathbb{S}^2} |(v-u)\cdot \omega|\sqrt{\mu(v)\mu(u)}f(u)\,d\omega du,
\\
(K_2f)(v)&=\int_{\mathbb{R}^3}\int_{\mathbb{S}^2} |(v-u)\cdot \omega|\sqrt{\mu(u)\mu(u')}f(v')\,d\omega du,\nonumber\\
&\quad+\int_{\mathbb{R}^3}\int_{\mathbb{S}^2}B(v-u,\t)\sqrt{\mu(u)\mu(v')}f(u')\,d\omega du,
\notag
\\
\nu(v)&=\int_{\mathbb{R}^3}\int_{\mathbb{S}^2}B(v-u,\t)\mu(u)\,d\omega du\cong 1+|v|.
\notag
\end{align}

We list some properties of $\FL, \Gamma$ for later use:
\begin{enumerate}[(I)]
\item The null space of the operator $\FL$  is the 5-dimensional space of collision invariants:
$$
 \mathbb{N}=\mbox{Ker}\,\FL=span \{\sqrt{\mu},\quad  (v-\fu)\sqrt{\mu},\quad (|v-\fu|^2-3)\sqrt{\mu}\}.
$$
And let $\bf{P}$ denote the projection operator from $L^2(\mathbb{R}^3)$ to $ \mathbb{N}$.

\smallskip

\item $\nu(v)$ satisfies
\begin{equation*}
 \nu_0(1+|v|)\leq \nu(v)\leq  \nu_1 (1+|v|).
\end{equation*}

\smallskip

\item The operator $K$ satisfies the following Grad's estimates
	\begin{align*}
	Kf(v)=\int_{\mathbb{R}^3}k(v,\eta)f(\eta)\,d\eta,\quad
	\end{align*}
	where  $k(v,\eta)$ 
	\begin{align}\label{2.15}
	0\leq |k(v,\eta)|\leq \f{C}{|v-\eta|}e^{-\f{|v-\eta|^2}{8}}e^{-\f{||v|^2-|\eta|^2|^2}{8|v-\eta|^2}}+C|v-\eta| e^{-\f{|v|^2}{4}}e^{-\f{|\eta|^2}{4}},
	\end{align}
	where $C>0$ is a given constant.
Following \eqref{2.15}, it is direct to have
\begin{align*}
\int_{\mathbb{R}^3}\Big|k(v,\eta)\cdot\f{(1+|v|)^{\a}}{(1+|\eta|)^{\a}}\Big|d\eta&\leq  C_{\alpha}(1+|v|)^{-1}.
\end{align*}

\smallskip

\item $\FL$ satisfies
\begin{equation*}
\begin{split}
\int_{\mathbb{R}^3}g\FL g dv &\geq c_0\|(\mathbf{I-P}) g \|_{\nu}^2,\\
\int_{\mathbb{R}^3}\nu |\FL^{-1} h|^2 dv &\leq \bar{c}_0 \int_{\mathbb{R}^3}\nu^{-1} | h|^2 dv,
\end{split}
\end{equation*} for $h\in \mathbb{N}^{\perp}$.

\smallskip

\item The nonlinear term $\Gamma(f, f)\in \mathbb{N}^{\perp}$
and
\begin{align*}
\|\nu^{-1} w \Gamma(f,f)\|_{L^\infty_v}\leq C \|w f \|^2_{L^\infty_v}.
\end{align*}
\end{enumerate}

\

We introduce a lemma which will be used to obtain the uniform $L^\infty_{x,v}$ of approximate solutions.
\begin{lemma}\cite{DHWZ} \label{lemA.1}
Consider  a sequence $\{a_i\}_{i=0}^\infty $  with each $a_i\geq0$. 
For any fixed $k\in\mathbb{N}_+$, we denote $$A_i^k=\max\{a_i, a_{i+1},\cdots, a_{i+k}\}.$$
	
\noindent{(1)} Assume $D\geq0$.  If $a_{i+1+k}\leq \f18 A_i^{k}+D$ for $i=0,1,\cdots$, then it holds that
\begin{equation*}
A_i^k\leq \left(\f18\right)^{\left[\frac{i}{k+1}\right]}\cdot\max\{A_0^k, \ A_1^k, \cdots, \ A_k^k \}+\f{8+k}{7} D,\quad\mbox{for}\quad i\geq k+1.
\end{equation*}
	
\noindent{(2)} 
Let $0\leq \eta<1$ with $\eta^{k+1}\geq\frac14$.  If $a_{i+1+k}\leq \f18 A_i^{k}+C_k \cdot \eta^{i+k+1}$ for $i=0,1,\cdots$, then it holds that
\begin{align*}
A_i^k\leq \left(\f18\right)^{\left[\frac{i}{k+1}\right]}\cdot\max\{A_0^k, \ A_1^k, \cdots, \ A_k^k \}+2C_k\f{8+k}{7} \eta^{i+k},\quad\mbox{for}\quad i\geq k+1.
\end{align*}
\end{lemma}

\


\section{Existence for the Linearized problem }\label{L}

This section is devoted to the existence result for the following linearized problem with a source term
\begin{align}\label{1.7-1}
\begin{cases}
v_3 \partial_x \mathfrak{f}+\FL \mathfrak{f}=\tilde{S},\\
\mathfrak{f}(0,v)|_{v_3>0}=\mathfrak{f}(0,Rv)+f_b(Rv),\\
\lim_{x\rightarrow\infty}\mathfrak{f}(x,v)=0,
\end{cases}
(x,v)\in \R_+\times \R^3.
\end{align}
The main result is
\begin{theorem}\label{thm3.1}
Recall the weight function $w(v)$ in \eqref{wt}, and let $\beta\geq 3$ and $0\leq \varsigma<\f14$. We assume   \eqref{1.12} and
\begin{align}\label{3.2}
\tilde{S}\in \mathbb{N}^{\perp}, \quad |wf_b|_{L^\infty_{v}}+\|\nu^{-1}w(v) e^{\sigma_0x}\tilde{S}\|_{L^\infty_{x,v}}<\infty,
\end{align}
then there exists a  unique solution $f$ to \eqref{1.7-1} such that
\begin{align}\label{3.132-0}
\|e^{\s x} w \mathfrak{f}\|_{L^\infty_{x,v}} +|w\mathfrak{f}(0,\cdot)|_{L^\infty_{v}}
\leq \frac{\tilde{C}}{\sigma_0-\sigma}\{\|e^{\s_0x} \nu^{-1}w\tilde{S}\|_{L^\infty_{x,v}}+|wf_b|_{L^\infty_v}\},
\end{align}
and
\begin{align}\label{3.132-01}
\|e^{\s x}\mathfrak{f}\|_{L^2_{x,v}}
\leq \frac{\tilde{C}}{\sigma_0-\sigma}\{\|e^{\s_0x}\tilde{S}\|_{L^2_{x,v}}+|(1+|v|)f_b|_{L^2_v}\},
\end{align}
where $\tilde{C}>0$ is a positive constant independent of $\sigma\in(0,\sigma_0)$. Moreover, if
$\tilde{S}$ is continuous in $\R_+\times \mathbb{R}^3$ and $f_b(v)$ is continuous in $\{v\in\R^3\}$, then $f(x,v)$ is continuous away from the  grazing set $\big\{(0,v)\ : \ v\in\R^3,\  v_3=0\big\}$.
\end{theorem}
Let $\chi(x)$ be a monotonic  smooth cut-off function
\begin{equation*}
\chi(x)\equiv 1,\  \mbox{for} \ x\in[0,1], \quad\mbox{and}\quad \chi(x)\equiv 0,\  \mbox{for} \  x\in[2,+\infty).
\end{equation*}
Similar as in \cite{GPS-1988}, we define
\begin{equation*}
f(x,v):=\mathfrak{f}(x,v)+ \chi(x) f_b(v),
\end{equation*}
then \eqref{1.7-1} is equivalent
\begin{align}\label{1.7-2}
\begin{cases}
v_3 \partial_x f+\FL f=g:=\tilde{S}+v_3\partial_x\chi(x)f_b(v)+\chi(x) \FL f_b,\\
f(0,v)|_{v_3>0}=f(0,Rv),\\
\lim_{x\rightarrow\infty} f(x,v)=0.
\end{cases}
\end{align}

To prove the existence of solution to \eqref{1.7-2}, we first need to consider a truncated problem. We denote $\Omega_d:=(0,d)$ with $d\geq1$ and denote the phase boundary of  $\Omega_d \times \mathbb{R}^{3}$ as
$\gamma =\partial \Omega_d \times \mathbb{R}^{3}$. We split $\gamma$ into three disjoint parts, outgoing
boundary $\gamma _{+},$ the incoming boundary $\gamma _{-},$ and the
singular boundary $\gamma _{0}$ for grazing velocities$:$
\begin{align}
\gamma _{+} &=\{(x,v)\in \partial \Omega_d \times \mathbb{R}^{3}:
n(x)\cdot v>0\}, \nonumber\\
\gamma _{-} &=\{(x,v)\in \partial \Omega_d \times \mathbb{R}^{3}:
n(x)\cdot v<0\}, \nonumber\\
\gamma _{0} &=\{(x,v)\in \partial \Omega_d \times \mathbb{R}^{3}:
n(x)\cdot v=0\}.\nonumber
\end{align}
where $\vec{n}(x)$ is the outward unit normal. It is direct to know that $\partial\Omega_d=\{0,d\}$, $\vec{n}(0)=(0,0,-1)$ and $\vec{n}(d)=(0,0,1)$.

Now we consider the truncated   problem with penalized term
\begin{align}\label{S3.3}
\begin{cases}
\v f^\v+ v_3 \partial_x f^\v+\FL f^\v=g,\\
f^\v(x,v)|_{\gamma_-}=f^\v(x,R_xv),\\
\end{cases}
(x,v)\in \Omega_d\times \R^3,
\end{align}
where $\v\in(0,1]$ and $R_xv:=v-2(v\cdot \vec{n}(x)) \vec{n}(x)$. We also define
\begin{equation*}
h^\v(x,v):=w(v) f^\v(x,v), 
\end{equation*}
then 
\eqref{S3.3} can be rewritten as
\begin{align}\label{S3.8}
\begin{cases}
\dis \v h^\v+v_3\partial_x h^\v+\nu(v) h^\v=K_{w} h^\v+wg,\\[2mm]
\dis h^\v(x,v)|_{\gamma_-}=h^\v(x,R_xv),
\end{cases}
\end{align}
where $K_wh=wK(\f{h}{w}).$ It is direct to know that
\begin{align}\label{3.8}
K_wh(v)=\int_{\R^3} k_w(v,u) h(u) du\quad\mbox{with}\quad
k_w(v,u)=w(v) k(v,u) w(u)^{-1}.
\end{align}


\subsection{A priori $L^\infty_{x,v}$ estimate}
For the approximate problem  \eqref{S3.8},  the most difficult part is to obtain the $L^\infty_{x,v}$-bound.

\begin{definition}
Given $(t,x,v),$ let $[X(s),V(s)]$ 
be the backward characteristics for \eqref{S3.8}, which is determined by
\begin{align*}
\begin{cases}
\dis \frac{{\dd}X(s)}{{\dd}s}=V_3(s),
\quad \frac{{\dd}V(s)}{{\dd}s}=0,\\[2mm]
[X(t),V(t)]=[x,v].
\end{cases}
\end{align*}
The solution is then given by
\begin{equation*}
[X(s),V(s)]=[X(s;t,x,v),V(s;t,x,v)]=[x-(t-s)v_3,v].
\end{equation*}
\end{definition}
Now for each $(x,v)$ with $x\in \bar{\Omega}_d$ and $v_3\neq 0,$ we define its {backward exit time} $t_{\mathbf{b}}(x,v)\geq 0$ to be the last moment at which the
back-time straight line $[X({-\tau};0,x,v),V({-\tau};0,x,v)]$ remains in $\bar{\Omega}$:
\begin{equation*}
t_{\mathbf{b}}(x,v)={\sup\{s \geq 0:x-\tau v_3\in\bar{\Omega}_d\text{ for }0\leq \tau\leq s\}.}
\end{equation*}
We also define
\begin{equation*}
x_{\mathbf{b}}(x,v)=x(t_{\mathbf{b}})=x-t_{\mathbf{b}}(x,v)\, v_3\in \partial \Omega_d .
\end{equation*}
We point out that $X(s)$, $t_{\mathbf{b}}(x,v)$ and $x_{\mathbf{b}}(x,v)$ are independent of  the horizontal velocity $v_h:=(v_1,v_2)$.

Let $x\in \bar{\Omega}_d$, $(x,v)\notin \gamma _{0}\cup \g_{-}$ and
$
(t_{0},x_{0},v_{0})=(t,x,v)$, and inductively define for $k\geq1$
\begin{equation*}
(t_{k+1},x_{k+1},v_{k+1})=(t_{k}-t_{\mathbf{b}}(x_{k},v_{k}),x_{\mathbf{b}}(x_{k},v_{k}), R_{x_{k+1}}v_{k}).
\end{equation*}
We define the back-time cycle as
\begin{equation}\label{2.7}
\left\{\begin{aligned}
X_{cl}(s;t,x,v)&=\sum_{k}\Fi_{[t_{k+1},t_{k})}(s)\{x_{k}-v_{k,3}\cdot(t_{k}-s)\},\\[1.5mm]
V_{cl}(s;t,x,v)&=\sum_{k}\Fi_{[t_{k+1},t_{k})}(s)v_{k},
\end{aligned}\right.
\end{equation}
Clearly, for $k\geq 1$ and $(x,v)\notin \g_0\cup \g_-$, it holds that
\begin{equation}\label{2.8}
\begin{split}
&x_k=\frac{1-(-1)^k}{2} x_1+\frac{1+(-1)^k}{2} x_2,\quad  v_{k,h}=v_{0,h}, \quad v_{k,3}=(-1)^{k} v_{0,3},\\
&t_k-t_{k+1}=t_1-t_2=\frac{d}{|v_{0,3}|}>0,\quad \nu(v)\equiv \nu(v_{k}).
\end{split}
\end{equation}

We can represent the solution of \eqref{S3.8} in a mild formulation which enables us to get the $L^\infty$ bound of solutions.  Indeed, for later use, we consider the following iterative linear problems involving a parameter $\la\in [0,1]$:
\begin{equation}\label{2.9}
\begin{cases}
\v h^{i+1}+v_3\partial_x h^{i+1}+\nu(v) h^{i+1}=\lambda K_w h^i +wg,\\[3mm]
\dis h^{i+1}(x,v)|_{\gamma_-}= h^i(x,R_xv)+w(v)r(x,v),
\end{cases}
\end{equation}
for $i=0,1, 2,\cdots$, where $h^0\equiv0$ and $w(v) r(x,v)\in L^\infty(\gamma_-)$ is any given function . For the mild formulation of \eqref{2.9}, we have the following lemma whose proof is omitted for brevity as it is similar to that in \cite{Guo2}.

\begin{lemma}
Let $0\leq \lambda\leq 1$. For each $(x,v)\in \bar{\Omega}_d\times \mathbb{R}^3\setminus (\gamma_0\cup \gamma_-)$,
we have
\begin{align}\label{2.10}
h^{i+1}(x,v)&=e^{-\nu_\v(v) (t-t_{k})} h^{i-k+1}(x_k,v_k)+\sum_{l=0}^{k-1} e^{-\nu_\v(v) (t-t_{l+1})} wr(x_{l+1},v_l)  \nonumber\\
&\quad+\sum_{l=0}^{k-1} \int_{t_{l+1}}^{t_l}e^{-\nu_\v(v) (t-s)} [\lambda K_w h^{i-l}+wg](X_{cl}(s),v_{l}) ds,
\end{align}
where and whereafter denote $\nu_\v(v)=\v+\nu(v)$, and  the parameter $k\gg 1$ is the collision times of the particle with  boundary $\partial\Omega_d$.
\end{lemma}


\begin{lemma}\label{lemS3.3}
Let $h^i$, $i=0,1,2,\cdots$, be the solutions to \eqref{2.9}, satisfying
$$
\|h^i\|_{L^\infty_{x,v}}+|h^i|_{L^\infty(\gamma_+)}<\infty.
$$ 
Then there exists $k_0>0$ large enough  such that for $i\geq 2k_0$, it holds 
that
\begin{align}\label{3.13}
\|h^{i+1}\|_{L^\infty_{x,v}}+|h^{i+1}|_{L^\infty(\g_+)}&\leq \frac18 \sup_{0\leq l\leq 2k_0} \{\|h^{i-l}\|_{L^\infty_{x,v}}+|h^{i-l}|_{L^\infty(\g_+)}\}\nonumber\\
&\hspace{-3mm} +C\{\|\nu^{-1}wg\|_{L^\infty_{x,v}}+|wr|_{L^\infty(\gamma_-)}\}+C \sup_{0\leq l\leq 2k_0}\left\{\left\|\frac{h^{i-l}}{w}\right\|_{L^2_{x,v}}\right\}.
\end{align}
Moreover, if $h^i\equiv h$ for $i=1,2,\cdots$, i.e., $h$ is a solution, then \eqref{3.13} is reduced to the following estimate
\begin{align}\label{3.14}
\|h\|_{L^\infty_{x,v}}{+|h|_{L^\infty(\g_+)}}&\leq C(\|\nu^{-1}wg\|_{L^\infty_{x,v}}+|wr|_{L^\infty(\gamma_-)})+C\left\|\frac{h}{w}\right\|_{L^2_{x,v}}.
\end{align}
We emphasize that the positive constant  $C>0$  depends on $k_0$, and is  independent of $d$, $\lambda\in[0,1]$ and $\vep>0$.
\end{lemma}

\noindent{\bf Proof.}
For $(x,v)\notin \g_0\cup \g_-$, it is noted that
\begin{equation*}
\frac{\nu_{\v}(v)}{|v_{0,3}|}=\frac{\v+\nu(v)}{|v_{0,3}|}\geq \nu_0 \frac{1+|v|}{|v|}\geq \nu_0>0,
\end{equation*}
which, together with \eqref{2.8}, yields that
\begin{align*}
\left|e^{-\nu_\v(v) (t-t_{k})} h^{i-k+1}(x_k,v_k)\right|
&\leq \left|e^{-\nu_\v(v) (t_1-t_{k})} h^{i-k+1}(x_k,v_k)\right|\nonumber\\
&\leq C |h^{i-k+1}|_{L^\infty(\gamma_+)}  \exp\left\{-\nu_{\v}(v)(k-1)\frac{d }{|v_{0,3}|}\right\}\nonumber\\
&\leq C |h^{i-k+1}|_{L^\infty(\gamma_+)}  \exp\left\{ -\nu_0 d (k-1) \right\}\nonumber\\
&\leq C |h^{i-k+1}|_{L^\infty(\gamma_+)}  e^{ -\frac12 \nu_0 d k }.
\end{align*}
For the second term on RHS of \eqref{2.10}, one has that
\begin{align*}
\left|\sum_{l=0}^{k-1} e^{-\nu_\v(v) (t-t_{l+1})} wr(x_{l+1},v_l) \right|\leq Ck |wr|_{L^\infty(\gamma_-)}.
\end{align*}
The last term on RHS of \eqref{2.10} is bounded by
\begin{align*}
\left|\sum_{l=0}^{k-1} \int_{t_{l+1}}^{t_l}e^{-\nu_\v(v) (t-s)}  wg(X_{cl}(s),v_{l}) ds\right|
&\leq C\|\nu^{-1}wg\|_{L^\infty_{x,v}} \sum_{l=0}^{k-1} \int_{t_{l+1}}^{t_l}e^{-\nu_\v(v) (t-s)} \nu(v)ds\nonumber\\
&\leq C\|\nu^{-1}wg\|_{L^\infty_{x,v}}.
\end{align*}

For the third term on RHS of \eqref{2.10}, we use \eqref{2.10} again to obtain
\begin{align}\label{2.17}
&\sum_{l=0}^{k-1} \int_{t_{l+1}}^{t_l}e^{-\nu_\v(v) (t-s)}  \lambda K_w h^{i-l}(X_{cl}(s),v_{l,3}) ds\nonumber\\
&=\lambda\sum_{l=0}^{k-1} \int_{t_{l+1}}^{t_l}e^{-\nu_\v(v) (t-s)}  \int_{\R^3} k_w(v_l,v') h^{i-l}(X_{cl}(s),v') dv'ds\nonumber\\
&\leq \lambda^2 \sum_{l=0}^{k-1} \int_{t_{l+1}}^{t_l}e^{-\nu_\v(v) (t-s)}  \int_{\R^3} |k_w(v_l,v')| dv'ds\nonumber\\
&\qquad\times \sum_{j=0}^{k-1} \int_{t'_{j+1}}^{t'_j}e^{-\nu_\v(v') (s-s_1)}\int_{\R^3} |k_w(v'_j,v'') h^{i-l-j-1}(X_{cl}'(s_1),v'')|dv'' ds_1 \nonumber\\
&\quad+  C\Big( e^{-\frac12\nu_0 dk} \sup_{0\leq l\leq 2k}|h^{i-l}|_{L^\infty(\gamma_+)}+\|\nu^{-1}wg\|_{L^\infty_{x,v}}+k|wr|_{L^\infty(\gamma_-)}\Big),
\end{align}
where we have denoted $X_{cl}'(s_1)=X_{cl}(s_1;s,X_{cl}(s),v')$, and $t'_j,v'_j$ are the corresponding times and velocities for specular cycles.

For the first term on RHS of \eqref{2.17}, we divide the proof into several cases:

\noindent{\it Case 1.} For $|v|\geq N$, the first term on RHS of \eqref{2.17} is bounded by
\begin{align*}
&\sum_{l=0}^{k-1} \int_{t_{l+1}}^{t_l}e^{-\nu_\v(v) (t-s)}  \int_{\R^3} |k_w(v_l,v')| dv'ds \cdot \sup_{0\leq l\leq 2k} \|h^{i-l}\|_{L^\infty_{x,v}}\nonumber\\
&\leq C\sup_{0\leq l\leq 2k} \|h^{i-l}\|_{L^\infty_{x,v}}\sum_{l=0}^{k-1} \int_{t_{l+1}}^{t_l}e^{-\nu_\v(v) (t-s)} \frac{1}{1+|v_l|} ds\nonumber\\
&\leq C\frac1{1+|v|} \sup_{0\leq l\leq 2k} \|h^{i-l}\|_{L^\infty_{x,v}}
\leq \frac{C}{N}\sup_{0\leq l\leq 2k} \|h^{i-l}\|_{L^\infty_{x,v}},
\end{align*}
where we have used the fact $|v|\equiv |v_l|$ for $l=0,1,\cdots$.

\

\noindent{\it Case 2.} For either $|v|\leq N, |v'|\geq 2N$ or $|v'|\leq 2N, |v''|\geq 3N$,  noting $|v_l|=|v|$ and $|v'_j|=|v'|$ we get either $|v_l-v'|\geq N$ or $|v'_j-v''|\geq N$, then either one of the following is valid for some small positive constant $0<c_1\leq \frac{1}{32}$:
\begin{equation}\label{3.21}
\begin{split}
\dis |k_w(v_l,v')|&\leq e^{-c_1N^{2}}|k_w(v_l,,v') \exp{\left(  c_1|v_l-v'|^{2}\right)}|,\\[1mm]
\dis |k_w(v'_j,v'')|&\leq e^{-c_1N^{2}}|k_w(v_j',,v'')\exp{\left(  c_1|v'_j-v'|^{2}\right)}|,
\end{split}
\end{equation}
which, together with \eqref{2.15},\eqref{3.8}, yields that
\begin{equation}\label{3.22}
\begin{split}
\dis \int_{\R^3}|k_w(v_l,v')e^{c_1|v_l-v'|^{2}}|dv'&\leq \frac{C}{1+|v|},\\
\dis \int_{\R^3}|k_w(v'_j,v'')e^{c_1|v_j'-v''|^{2}}|dv''&\leq \frac{C}{1+|v'|}.
\end{split}
\end{equation}
Using \eqref{3.21}-\eqref{3.22}, one has
\begin{align*}
&\sum_{l=0}^{k-1} \int_{t_{l+1}}^{t_l}e^{-\nu_\v(v) (t-s)}  \left\{\iint_{|v|\leq N, |v'|\geq 2N} + \iint_{|v'|\leq 2N, |v''|\geq 3N}  \right\} (\cdots)dv'' ds_1 dv'ds\nonumber\\
&\leq Ce^{-c_1 N^2} \sup_{0\leq l\leq 2k} \|h^{i-l}\|_{L^\infty_{x,v}}\leq  \frac{C}{N}\sup_{0\leq l\leq 2k} \|h^{i-l}\|_{L^\infty_{x,v}},
\end{align*}

\noindent{\it Case 3.} For either $|v|\leq N, |v'|\leq 2N$, $|v''|\leq 3N$, this is the last remaining case. We denote $D=\{|v'|\leq 2N,|v''|\leq 3N\}$. Noting $\nu_{\v}(v)\geq \nu_{0}$,  the corresponding part is bounded by
\begin{align}
&\sum_{l=0}^{k-1} \int_{t_{l+1}}^{t_l}e^{-\nu_0 (t-s)}  \iint_{D} |k_w(v_l,v')k_w(v'_j,v'')| dv'' dv'ds\nonumber\\
&\qquad\times \sum_{j=0}^{k-1} \left(\int_{t'_{j}-\frac1N}^{t'_j}+\int_{t'_{j+1}}^{t'_j-\frac1N}\right)e^{-\nu_0 (s-s_1)} | h^{i-l-j-1}(X_{cl}'(s_1),v'')| ds_1 \nonumber\\
&\leq \sum_{l=0}^{k-1} \int_{t_{l+1}}^{t_l}e^{-\nu_0 (t-s)}  \iint_{D} |k_w(v_l,v')k_w(v'_j,v'')| dv'' dv'ds\nonumber\\
&\qquad\times \sum_{j=0}^{k-1} \int_{t'_{j+1}}^{t'_j-\frac1N} e^{-\nu_0 (s-s_1)} | h^{i-l-j-1}(X_{cl}'(s_1),v'')| ds_1 +C\frac{k}{N}\sup_{0\leq l\leq 2k} \|h^{i-l}\|_{L^\infty_{x,v}} \nonumber\\
&\leq \sum_{l=0}^{k-1} \int_{t_{l+1}}^{t_l}e^{-\nu_0 (t-s)} ds \left(\iint_{D} \sum_{j=0}^{k-1} \int_{t'_{j+1}}^{t'_j-\frac1N} e^{-\nu_0 (s-s_1)}  |k_w(v_l,v')k_w(v'_j,v'')|^2  ds_1dv'dv''\right)^{\frac12}\nonumber\\
&\qquad\times \left( \iint_{D} \sum_{j=0}^{k-1}\int_{t'_{j+1}}^{t'_j-\frac1N} e^{-\nu_0 (s-s_1)} | h^{i-l-j-1}(x'_j-v'_{j,3} (t'_j-s_1),v'')|^2 ds_1dv'dv''\right)^{\frac12} \nonumber\\
&\qquad+C\frac{k}{N}\sup_{0\leq l\leq 2k} \|h^{i-l}\|_{L^\infty_{x,v}}.\label{3.24}
\end{align}
It follows from \eqref{2.15}  that
\begin{align}\nonumber
\iint_{D} \sum_{j=0}^{k-1} \int_{t'_{j+1}}^{t'_j-\frac1N} e^{-\nu_0 (s-s_1)}  |k_w(v_l,v')k_w(v'_j,v'')|^2  ds_1dv'dv''\leq C.
\end{align}
Define $y:=x'_j-v'_{j,3} (t'_j-s_1)$. We have  $x'_j=0\,  \mbox{or}\, d$ and $v'_{j,3}=(-1)^{j} v'_{0,3}$. For $t'_j=t'_j(s_1;s_1, X_{cl}(s),v')$, it holds that
\begin{align}\nonumber
s-t'_j=
\begin{cases}
\displaystyle \frac{X_{cl}(s)}{|v'_{0,3}|}+(j-1) \frac{d}{|v'_{0,3}|},\quad \mbox{for}\  v'_{0,3}>0,\\[3mm]
\displaystyle \frac{d-X_{cl}(s)}{|v'_{0,3}|}+(j-1) \frac{d}{|v'_{0,3}|},\quad \mbox{for}\  v'_{0,3}<0,
\end{cases}
\end{align}
which yields that
\begin{align}\nonumber
y=
\begin{cases}
x'_j-(-1)^{j}\Big\{v'_{0,3}(s-s_1)-[X_{cl}(s)+(j-1)d] \Big\},\quad  \mbox{for}\  v'_{0,3}>0,\\[2mm]
x'_j-(-1)^{j} \Big\{v'_{0,3}(s-s_1)+[j d-X_{cl}(s)]\Big\},\quad  \mbox{for}\  v'_{0,3}>0.
\end{cases}
\end{align}
Since  $x'_j=0\,  \mbox{or}\, d$, which is independent of $v_{0,3}'$, thus we have
\begin{equation}\nonumber
\left|\frac{dy}{d v'_{0,3}}\right| = (s-s_1)\geq \frac1N, \quad \mbox{for} \  s_1\in [t'_{j+1},t'_j-\frac1N],
\end{equation}
which yields that
\begin{align}
& \left( \iint_{D} \sum_{j=0}^{k-1}\int_{t'_{j+1}}^{t'_j-\frac1N} e^{-\nu_0 (s-s_1)} | h^{i-l-j-1}(x'_j-v'_{j,3} (t'_j-s_1),v'')|^2 ds_1dv'dv''\right)^{\frac12} \nonumber\\
&\leq C_N k \sup_{0\leq l\leq 2k}\|\frac{h^{i-l}}{w}\|_{L^2_{x,v}}.\nonumber
\end{align}
Then the RHS of \eqref{3.24} is bounded by
\begin{equation*}
C\frac{k}{N}\sup_{0\leq l\leq 2k} \|h^{i-l}\|_{L^\infty_{x,v}}+C_{k,N} \sup_{0\leq l\leq 2k}\|\frac{h^{i-l}}{w}\|_{L^2_{x,v}}.
\end{equation*}

Combining above estimates, we obtain
\begin{align*}
\|h^{i+1}\|_{L^\infty_{x,v}}+|h^{i+1}|_{L^\infty(\gamma_+)}
&\leq C(e^{-\frac12\nu_0 d k}+\frac{k}{N}) \sup_{0\leq l\leq 2k} \{\|h^{i-l}\|_{L^\infty_{x,v}}+|h^{i-l}|_{L^\infty(\gamma_+)}\}\nonumber\\
&\, +C_{N,k} \sup_{0\leq l\leq 2k}\|\frac{h^{i-l}}{w}\|_{L^2_{x,v}}+C_{k} (\|\nu^{-1}wg\|_{L^\infty_{x,v}}+|wr|_{L^\infty(\gamma_-)}).
\end{align*}
First taking $k$ large , and then letting $N$ suitably large so that
\begin{equation}\nonumber
C(e^{-\frac12\nu_0 k}+\frac{k}{N})\leq \frac18,
\end{equation}
which implies \eqref{3.13}. This completes  the proof of Lemma \ref{lemS3.3}. $\hfill\Box$


\subsection{Approximate solutions and uniform estimate}
Now we are in a position to  construct solutions to \eqref{S3.3} or equivalently \eqref{S3.8}. First of all, we consider the following approximate problem
\begin{align}\label{S3.54}
\begin{cases}
\dis \v f^n+ v_3 \partial_xf^n+\nu(v) f^n-Kf^n=g,\\[2mm]
\dis f^n(x,v)|_{\gamma_{-}}=(1-\frac{1}{n}) f^n(x,R_xv),
\end{cases}
(x,v)\in \Omega_d\times \mathbb{R}^3,
\end{align}
where $\v\in(0,1]$ is arbitrary and $n>1$ is an integer. For later use, we choose $n_0>1$ large enough such that
$$
\frac18 (1-\frac1n)^{-\frac{k_0+1}{2}}\leq \frac12
$$
for any $n\geq n_0$, where  $k_0>0$ is the one fixed  in Lemma \ref{lemS3.3}.

\begin{lemma}\label{lemS3.4}
Let $\v>0, d\geq 1$, $n\geq n_0$, and $\beta\geq 3$. Assume $\|\nu^{-1}wg\|_{L^\infty_{x,v}}<\infty$.
Then there exists a unique solution $f^n$ to \eqref{S3.54} satisfying
\begin{align*}
\|wf^{n}\|_{L^\infty_{x,v}}+|wf^n|_{L^\infty(\gamma_+)}\leq C_{\v,n} \|\nu^{-1}wg\|_{L^\infty_{x,v}},
\end{align*}
where the positive constant $C_{\v,n}>0$ depends only on $\v$ and $n$. Moreover, if
g is continuous in $\Omega_d\times\mathbb{R}^3$, then $f^n$ is continuous away from grazing set $\gamma_0$.
\end{lemma}

\noindent{\bf Proof.}  We consider the solvability of the  following boundary value problem
\begin{equation}\label{S3.56-1}
\begin{cases}
\mathcal{L}_\lambda f:=\v f+v_3 \partial_x f+\nu(v) f-\lambda Kf=g,\\
f(x,v)|_{\gamma_-}=(1-\frac1n) f(x,R_xv),
\end{cases}
\end{equation}
for $\lambda\in[0,1]$. For brevity we denote $\mathcal{L}_\lambda^{-1}$ to be the solution operator associated with the problem, meaning that $f:=\mathcal{L}_\lambda^{-1} g$ is a solution to the BVP \eqref{S3.56-1}. Our idea is to prove the existence of $\mathcal{L}_0^{-1}$, and then extend to obtain the existence of  $\mathcal{L}_1^{-1}$ by a continuous argument on $\la$.  Since the proof is very long, we split it  into  several steps.   

\medskip
\noindent{\it Step 1.} In this step, we prove the existence of $\mathcal{L}_0^{-1}$. We consider the following approximate sequence
\begin{align}\label{S3.57}
\begin{cases}
\mathcal{L}_0 f^{i+1}=\v f^{i+1}+v_3\partial_x f^{i+1} + \nu(v) f^{i+1}=g,\\
f^{i+1}(x,v)|_{\gamma_-}=(1-\frac1n) f^{i}(x,R_xv),
\end{cases}
\end{align}
for $i=0,1,2,\cdots$, where we have set $f^0\equiv0$. We will construct $L^\infty$ solutions to \eqref{S3.57} for $i=0,1,2,\cdots$, and 
establish uniform $L^\infty$-estimates.

Firstly, we will solve inductively the linear equation \eqref{S3.57} by 
the method of characteristics. Let $h^{i+1}(x,v)=w(v)f^{i+1}(x,v)$.  For 
almost every $(x,v)\in\bar{\Omega}_d\times\mathbb{R}^3\backslash (\gamma_0\cup \gamma_-)$, 
one can write
\begin{align}\label{S3.58}
h^{i+1}(x,v)
&=e^{-(\v+\nu(v))t_{\mathbf{b}}}\cdot (1-\frac1n)  w(v) f^{i}(x_1,R_{x_1}v)\nonumber\\
&\quad+\int_{t_1}^t  e^{-(\v+\nu(v))(t-s)} (wg)(x-v_{0,3}(t-s),v) \dd s,
\end{align}
where $t_1=t-t_{\mathbf{b}}$.
We consider \eqref{S3.58} with $i=0$. Noting $h^0\equiv0$,   then 
it is straightforward to see that
\begin{align*}
\|h^{1}\|_{L^\infty_{x,v}}+|h^1|_{L^\infty(\gamma_+)}\leq   C \|\nu^{-1}wg\|_{L^\infty_{x,v}}<\infty.
\end{align*}
Therefore we have obtained the solution to \eqref{S3.57} with $i=0$.  Assume that we have already solved \eqref{S3.57} for $i\leq l$ and obtained
\begin{equation}\label{S3.62}
\|h^{l+1}\|_{L^\infty_{x,v}}+|h^{l+1}|_{L^\infty(\gamma_+)}\leq C_{l+1}  \|\nu^{-1}wg\|_{L^\infty_{x,v}}<\infty.
\end{equation}
We now consider \eqref{S3.57} for $i=l+1$. Noting \eqref{S3.62}, then we can solve \eqref{S3.57} by using \eqref{S3.58} with $i=l+1$. We still need to prove $h^{l+2}\in L^\infty$. Indeed, it follows from \eqref{S3.58} that
\begin{align*}
\|h^{l+2}\|_{L^\infty_{x,v}}+{|h^{l+2}|_{L^\infty(\g_+)}}&\leq C |h^{l+1}|_{L^\infty{(\g_+)}} +C\|\nu^{-1}wg\|_{L^\infty_{x,v}}\nonumber\\
&\leq C_{l+2}  \|\nu^{-1}wg\|_{L^\infty_{x,v}}<\infty.
\end{align*}
Therefore, inductively we have solved  \eqref{S3.57} for $i=0,1,2,\cdots$ and obtained
\begin{align}\label{S3.65}
\|h^{i}\|_{L^\infty_{x,v}}+|h^{i}|_{L^\infty(\gamma_+)}\leq C_{i} \|\nu^{-1}wg\|_{L^\infty_{x,v}}<\infty,
\end{align}
for $i=0,1,2,\cdots$.
The positive  constant $C_{i}$ may increase to infinity as $i\rightarrow \infty$.
Here, we emphasize that we first need to know the sequence $\{h^i\}_{i=0}^{\infty}$ is in $L^\infty_{x,v}$-space, otherwise one can not use Lemma \ref{lemS3.3} to get uniform $L^\infty_{x,v}$ estimates.

Since  $\Omega_d$ is a convex domain, let $(x,v)\in \Omega_d\times\mathbb{R}^3\backslash \gamma_0$, then it is easy to check that $t_{\mathbf{b}}(x,v)$ and $x_{\mathbf{b}}(x,v)$ are continuous. Therefore if $g$ and $r$ are continuous, we conclude that $f^i(x,v)$ is continuous away from grazing set.

\vspace{1mm}

Secondly, in order to take the limit $i\rightarrow\infty$, one has to get some uniform  estimates.  Multiplying \eqref{S3.57} by $f^{i+1}$ and integrating the resultant equality over $\Omega_d\times\mathbb{R}^3$, one obtains that
\begin{align}\label{S3.66}
& \v\|f^{i+1}\|^2_{L^2_{x,v}}+\frac12|f^{i+1}|^2_{L^2(\gamma_+)}+\|f^{i+1}\|^2_{\nu}\nonumber\\
&\leq \frac12(1-\frac1n)^2|f^{i}|^2_{L^2(\gamma_+)}+C\|g\|^2_{L^2_{x,v}}+ \frac{1}{4}\|f^{i+1}\|^2_{\nu}.
\end{align}
 Then, from \eqref{S3.66},  we have
\begin{align*}
2\v\|f^{i+1}\|^2_{L^2_{x,v}}+|f^{i+1}|^2_{L^2(\gamma_+)}+\|f^{i+1}\|^2_{\nu}
\leq (1-\frac1n)^2|f^{i}|^2_{L^2(\gamma_+)}+C \|g\|^2_{L^2_{x,v}}.
\end{align*}
Now we take the difference $f^{i+1}-f^i$ in \eqref{S3.57}, then by similar energy estimate as above, we obtain
\begin{align}\label{S3.68}
&2\v\|f^{i+1}-f^i\|^2_{L^2_{x,v}}+|f^{i+1}-f^i|^2_{L^2(\gamma_+)}+2\|f^{i+1}-f^i\|^2_{\nu}\nonumber\\
&\leq (1-\frac1n)^2|f^{i}-f^{i-1}|^2_{L^2(\gamma_+)} \leq \cdots \leq (1-\frac1n)^{2i} |f^1|^2_{L^2(\gamma_+)}\nonumber\\
&\leq  C(1-\frac1n)^{2i} \|g\|_{L^2_{x,v}}<\infty.
\end{align}
Noting $1-\frac1n<1$, thus $\{f^{i}\}_{i=0}^\infty$ is a Cauchy sequence in $L^2$, i.e.,
\begin{align*}
|f^{i}-f^j|^2_{L^2(\gamma_+)}+\|
{f^i-f^j}\|^2_{\nu}\rightarrow0,\quad\mbox{as} \  i,j\rightarrow\infty.
\end{align*}
And we also have, for $i=0,1,2,\cdots$, that
\begin{align}\label{S3.70}
|f^{i}|^2_{L^2(\gamma_+)}+\|f^{i}\|^2_{\nu}\leq C\|g\|^2_{L^2_{x,v}}.
\end{align}

\vspace{1mm}

Next we consider the uniform  $L^\infty_{x,v}$ estimate. Here we point out that Lemma \ref{lemS3.3} still holds by replacing 1 with $1-\frac1n$ in the boundary condition, and the constants in Lemma \ref{lemS3.3} do not depend on $n\geq1$.  Thus we apply Lemma \ref{lemS3.3} to obtain that
\begin{align*}
&\|h^{i+1}\|_{L^\infty_{x,v}}+|h^{i+1}|_{L^\infty{(\gamma_+)}}\nonumber\\
&\leq \frac18\sup_{0\leq l\leq 2k_0} \{\|h^{i-l}\|_{L^\infty_{x,v}}+|h^{i-l}|_{L^\infty(\gamma_+)}\}
+C\|\nu^{-1}wg\|_{L^\infty_{x,v}}+C\sup_{0\leq l\leq 2k_0}  \|f^{i-l}\|_{L^2_{x,v}} \nonumber\\
&\leq\frac18\sup_{0\leq l\leq 2k_0} \{\|h^{i-l}\|_{L^\infty_{x,v}}+|h^{i-l}|_{L^\infty(\gamma_+)}\}
+C_d\|\nu^{-1}wg\|_{L^\infty_{x,v}},
\end{align*}
where we have used \eqref{S3.70} in the second inequality.  Now we apply Lemma \ref{lemA.1} to obtain that for $i\geq 2k_0+1$,
\begin{align}\label{S3.72}
&\|h^{i}\|_{L^\infty_{x,v}}+|h^{i}|_{L^\infty{(\gamma_+)}}\nonumber\\
&\leq\left(\frac18\right)^{[\frac{i}{2k_0+1}]}  \max_{0\leq l\leq 2k_0} \Big\{\|h^{1}\|_{L^\infty_{x,v}}+|h^{1}|_{L^\infty{(\gamma_+)}}, \cdots, \|h^{2k_0}\|_{L^\infty_{x,v}}+|h^{2k_0}|_{L^\infty{(\gamma_+)}}\Big\}\nonumber\\
&\qquad+\frac{8+2k_0}{7} C_{d}\|\nu^{-1}wg\|_{L^\infty_{x,v}}\nonumber\\
&\leq C_{k_0,d}  \|\nu^{-1}wg\|_{L^\infty_{x,v}},
\end{align}
where we have used \eqref{S3.65} in the second inequality.
Hence it follows from \eqref{S3.72} and \eqref{S3.65} that
\begin{align}\label{S3.73}
\|h^{i}\|_{L^\infty_{x,v}}+|h^{i}|_{L^\infty{(\gamma_+)}}\leq C_{k_0,d}\|\nu^{-1}wg\|_{L^\infty_{x,v}}, \quad\mbox{for} \  i\geq 1.
\end{align}
Taking the difference $h^{i+1}-h^i$ and then applying Lemma \ref{lemS3.3} to $h^{i+1}-h^i$, we have that for $i\geq 2k_0$,
\begin{align}\label{S3.74}
&\|h^{i+2}-h^{i+1}\|_{L^\infty_{x,v}}+|h^{i+2}-h^{i+1}|_{L^\infty(\g_+)}\nonumber\\
&\leq \frac18 \max_{0\leq l\leq 2k_0} \Big\{\|h^{i+1-l}-h^{i-l}\|_{L^{\infty}_{x,v}}+|h^{i+1-l}-h^{i-l}|_{L^\infty(\g_+)}\Big\}\nonumber\\
&\qquad+C \sup_{0\leq l\leq 2k_0} \Big\{\|f^{i+1-l}-f^{i-l}\|_{L^2_{x,v}}\Big\}\nonumber\\
&\leq\frac18 \max_{0\leq l\leq 2k_0} \Big\{\|h^{i+1-l}-h^{i-l}\|_{L^{\infty}_{x,v}}+|h^{i+1-l}-h^{i-l}|_{L^\infty(\g_+)}\Big\}+C_{k_0,d}\|\nu^{-1}wg\|_{L^\infty_{x,v}} \eta_n^{i-2k_0}\nonumber\\
&\leq \frac18 \max_{0\leq l\leq 2k_0} \Big\{\|h^{i+1-l}-h^{i-l}\|_{L^{\infty}_{x,v}}+|h^{i+1-l}-h^{i-l}|_{L^\infty(\g_+)}\Big\}+C_{k_0,d} \|\nu^{-1}wg\|_{L^\infty_{x,v}}  \eta_n^{i+1},
\end{align}
where we have used \eqref{S3.68} and  denoted  $\eta_n:=
1-\frac1n
<1$. Here we choose $n$ large enough so that $\frac18 \eta_n^{-2k_0-1}\leq \frac12$, then it follows from \eqref{S3.74} and Lemma \ref{lemA.1}  that
\begin{align}\label{S3.75}
&\|h^{i+2}-h^{i+1}\|_{L^\infty_{x,v}}+|h^{i+2}-h^{i+1}|_{L^\infty(\gamma_+)}\nonumber\\
&\leq \left(\frac18\right)^{\left[\frac{i}{2k_0+1}\right]} \max_{0\leq l\leq 4k_0}\Big\{\|h^1\|_{L^\infty_{x,v}}+|h^{1}|_{L^\infty(\gamma_+)}, \cdots , \|h^{4k_0+1}\|_{L^\infty_{x,v}}++|h^{4k_0+1}|_{L^\infty(\gamma_+)}\Big\}\nonumber\\
&\quad+C_{k_0,d} \|\nu^{-1}wg\|_{L^\infty_{x,v}}\cdot  \eta_n^{i}\nonumber\\
&\leq C_{k_0,d} \|\nu^{-1}wg\|_{L^\infty_{x,v}}\cdot\Bigg\{\left(\frac18\right)^{\left[\frac{i}{2k_0+1}\right]}  +\eta_n^{i}\Bigg\},
\end{align}
for $i\geq 2k_0+1$. Then \eqref{S3.75} implies immediately that $\{h^i\}_{i=0}^\infty$ is a Cauchy sequence in $L^\infty_{x,v}$, i.e., there exists a limit function $
{h}\in L^\infty_{x,v}$ so that $\|h^{i}-{h}\|_{L^\infty_{x,v}}+|h^{i}-{h}|_{L^\infty(\gamma_+)}\rightarrow 0$ as $i\rightarrow \infty$. Thus we obtained a function
{$f:=\f{h}{w}$} solves
\begin{equation*}
\begin{cases}
\mathcal{L}_0 f=\v f+ v_3\partial_xf+\nu(v) f=g,\\[2mm]
\dis f(x,v)|_{\gamma_{-}}=(1-\frac{1}{n}) f(x,R_xv),
\end{cases}
\end{equation*}
with $n\geq n_0$ large enough.
Moreover, from \eqref{S3.73}, 
there exists a constant $C_{k_0,d}$ such that
\begin{align*}
\|h\|_{L^\infty_{x,v}}+|h|_{L^\infty(\g)}\leq C_{k_0,d} \|\nu^{-1}wg\|_{L^\infty_{x,v}}.
\end{align*}

\medskip
\noindent{\it Step 2.} {\it A priori estimates.} For any given $\lambda\in[0,1]$, let $f^n$ be the solution of \eqref{S3.56-1}, i.e.,
\begin{align}\label{S3.79}
\begin{cases}
\mathcal{L}_\lambda f^n=\v f^{n}+v_3\partial_x f^{n}+\nu(v) f^{n}-\lambda Kf^n=g,\\[2mm]
f^{n}(x,v)|_{\gamma_-}=(1-\frac1n) f^{n}(x,R_xv).
\end{cases}
\end{align}
Moreover we also assume that $\|wf^{n}\|_{L^\infty_{x,v}}+|wf^n|_{L^\infty(\gamma)}<\infty$. Firstly,  we shall consider   {\it a priori}  $L^2$-estimates. Multiplying \eqref{S3.79} by $f^{n}$, one has that
\begin{align}\label{S3.80}
& \v\|f^{n}\|^2_{L^2_{x,v}}+\frac12|f^{n}|^2_{L^2(\gamma_+)}-\frac12|f^{n}|^2_{L^2(\gamma_-)}+\|f^{n}\|^2_{\nu}\nonumber\\
&\leq \lambda\langle Kf^n, f^{n}\rangle +\frac{\v}{4} \|f^n\|_{L^2_{x,v}}+\frac{C}{\v}\|g\|^2_{L^2_{x,v}}.
\end{align}
We note that $\langle \FL f^n, f^{n}\rangle\geq0 $, which implies that
\begin{equation}\label{S3.81}
\langle Kf^n, f^{n}\rangle \leq \|f^n\|_{\nu}.
\end{equation}
On the other hand, it follows from $\eqref{S3.79}_2$ that
\begin{equation}\label{S3.82}
\frac12|f^{n}|^2_{L^2(\gamma_+)}-\frac12|f^{n}|^2_{L^2(\gamma_-)}= \frac12 |f^{n}|^2_{L^2(\gamma_+)}[1-(1-\frac1n)^2] \geq 0.
\end{equation}
Substituting \eqref{S3.81} and \eqref{S3.82} into \eqref{S3.80}, one has that
\begin{align}\label{S3.83}
\|\mathcal{L}^{-1}_\lambda g\|^2_{L^2_{x,v}}=\|f^{n}\|^2_{L^2_{x,v}}\leq  C_{\v} \|g\|^2_{L^2}.
\end{align}
Let $h^n:=w f^n$.  Then, by using \eqref{3.14} and \eqref{S3.83}, we  obtain
\begin{align}\label{S3.86}
\|w\mathcal{L}^{-1}_\lambda g\|_{L^\infty_{x,v}}+|w\mathcal{L}^{-1}_\lambda g|_{L^\infty(\gamma)}=\|h^n\|_{L^\infty_{x,v}}+|h^n|_{L^\infty(\gamma)}\leq C_{\v,k_0,d} \|\nu^{-1}wg\|_{L^\infty_{x,v}}.
\end{align}

\vspace{1mm}

On the other hand, Let $\nu^{-1}wg_1 \in L^\infty_{x,v}$ and $\nu^{-1}wg_2 \in L^\infty_{x,v}$.  Let  $f^n_1=\mathcal{L}^{-1}_\lambda g_1$ and $f^n_2=\mathcal{L}^{-1}_\lambda g_2$ be the solutions to \eqref{S3.79} with $g$ replaced by $g_1$ and $g_2$, respectively. Then we have that
\begin{align*}
\begin{cases}
\v (f^{n}_2-f^n_1)+v_3\partial_x (f^{n}_2-f^n_1)+\nu(v) (f^{n}_2-f^n_1)-\lambda K(f^{n}_2-f^n_1)=g_2-g_1,\\[2mm]
(f^{n}_2-f^n_1)(x,v)|_{\gamma_-}=(1-\frac1n) (f^{n}_2-f^n_1)(x,R_xv).
\end{cases}
\end{align*}
By similar arguments as in \eqref{S3.79}-\eqref{S3.86}, we obtain
\begin{equation}\label{S3.88}
\|\mathcal{L}^{-1}_\lambda g_2-\mathcal{L}^{-1}_\lambda g_1\|^2_{L^2_{x,v}}\leq  C_{\v,k_0,d} \|g_2-g_1\|^2_{L^2_{x,v}},
\end{equation}
and
\begin{equation}\label{S3.89}
\|w(\mathcal{L}^{-1}_\lambda g_2-\mathcal{L}^{-1}_\lambda g_1)\|_{L^\infty_{x,v}}
+|w(\mathcal{L}^{-1}_\lambda g_2-\mathcal{L}^{-1}_\lambda g_1)|_{L^\infty(\gamma)}\leq C_{\v,k_0,d}  \|\nu^{-1}w(g_2-g_1)\|_{L^\infty_{x,v}}.
\end{equation}
The uniqueness of solution to \eqref{S3.79} also follows from \eqref{S3.88}.
We point out that the constant $C_{\v,k_0,d}$ in \eqref{S3.83}, \eqref{S3.86}, \eqref{S3.88} and \eqref{S3.89} does not depend on $\lambda\in[0,1]$. This property is crucial for us to extend $\mathcal{L}_0^{-1}$ to $\mathcal{L}_1^{-1}$ by a bootstrap argument.

\medskip
\noindent{\it Step 3.}  In this step,  we shall prove the existence of solution $f^n$ to \eqref{S3.56-1} for sufficiently small $0<\lambda\ll1$, i.e., to prove the existence of operator $\mathcal{L}_{\lambda}^{-1}$. Firstly, we define the Banach space
\begin{equation}\nonumber
\mathbf{X}:=\Big\{f=f(x,v) :  \   wf\in L^\infty(\Omega\times\mathbb{R}^3), \ wf\in L^\infty(\gamma), \  \mbox{and} \  f(x,v)|_{\gamma_-}=(1-\frac1n) f(x,R_xv) \Big\}.
\end{equation}
Now we define
\begin{align}\nonumber
T_\lambda f=\mathcal{L}_0^{-1} \Big(\lambda K f+g\Big).
\end{align}
For any $f_1, f_2\in \mathbf{X}$, by using \eqref{S3.89}, we have  that
\begin{align}
&\|w(T_\lambda f_1-T_\lambda f_2)\|_{L^\infty_{x,v}}+|w(T_\lambda f_1-T_\lambda f_2)|_{L^\infty(\gamma)}\nonumber\\
&= \left\|w\{\mathcal{L}_0^{-1}(\lambda Kf_1+g)-\mathcal{L}_0^{-1}(\lambda Kf_2+g)\} \right\|_{L^\infty_{x,v}}\nonumber\\
&\qquad\quad+ \left|w\{\mathcal{L}_0^{-1}(\lambda Kf_1+g)-\mathcal{L}_0^{-1}(\lambda Kf_2+g)\} \right|_{L^\infty(\gamma)}\nonumber\\
&\leq C_{\v,k_0,d}\|\nu^{-1} w\{(\lambda Kf_1+g)-(\lambda Kf_2+g)\}\|_{L^\infty_{x,v}}\nonumber\\
&\leq \lambda C_{\v,k_0,d} \|w (Kf_1-Kf_2)\|_{L^\infty_{x,v}} \nonumber\\
&\leq \lambda C_{\v,k_0,d} \|w(f_1-f_2)\|_{L^\infty_{x,v}}.\nonumber
\end{align}
We take $\lambda_\ast>0$ sufficiently small such that $\lambda_\ast C_{\v,k_0,d}\leq 1/2$, then $T_\lambda : \mathbf{X}\rightarrow \mathbf{X}$ is a contraction mapping for $\lambda\in[0,\lambda_\ast]$. Thus $T_\lambda$ has a fixed point, i.e.,  $\exists\, f^\lambda\in \mathbf{X}$ such that
\begin{equation}\nonumber
f^\lambda=T_\lambda f^\lambda=\mathcal{L}_0^{-1} \Big(\lambda K f^\lambda+g\Big),
\end{equation}
which yields immediately that
\begin{align}\nonumber
\mathcal{L}_\lambda f^\lambda=\v f^\lambda + v_3 \partial_x f^\lambda+\nu f^\lambda-\lambda Kf^{\lambda}=g.
\end{align}
Hence, for any $\lambda\in[0,\lambda_\ast]$,  we have solved \eqref{S3.56-1} with $f^\lambda=\mathcal{L}_\lambda^{-1}g\in\mathbf{X}$.  Therefore we have obtained the existence of $\mathcal{L}_\lambda^{-1}$ for $\lambda\in[0,\lambda_\ast]$. Moreover the operator $\mathcal{L}_\lambda^{-1}$  has the properties \eqref{S3.83}, \eqref{S3.86}, \eqref{S3.88} and \eqref{S3.89}.

Next we define
\begin{align}\nonumber
T_{\lambda_\ast+\lambda}f=\mathcal{L}_{\lambda_\ast}^{-1}\Big(\lambda K f+g\Big).
\end{align}
Noting the estimates for $\mathcal{L}_{\lambda_\ast}^{-1}$ are independent of $\lambda_\ast$. By similar arguments, we can prove $T_{\lambda_\ast+\lambda} : \mathbf{X}\rightarrow \mathbf{X}$ is a contraction mapping  for $\lambda\in[0,\lambda_\ast]$.  Then we obtain the existence of operator $\mathcal{L}_{\lambda_\ast+\lambda}^{-1}$, and  \eqref{S3.83}, \eqref{S3.86}, \eqref{S3.88} and \eqref{S3.89}.   Step by step, we can finally obtain the existence of operator $\mathcal{L}_1^{-1}$, and $\mathcal{L}_1^{-1}$ satisfies the estimates in \eqref{S3.83}, \eqref{S3.86}, \eqref{S3.88} and \eqref{S3.89}. The continuity is easy to obtain since the convergence of sequence under consideration is always in $L^\infty_{x,v}$.  Therefore we complete the proof of Lemma \ref{lemS3.4}.  $\hfill\Box$

\begin{lemma}\label{lemS3.6}
Let $\v>0, d\geq 1$ and $\beta\geq3$, and  assume $\|\nu^{-1}wg\|_{L^\infty_{x,v}}<\infty$. Then there exists a unique solution $f^\v$ to solve the approximate linearized steady Boltzmann equation \eqref{S3.3}. Moreover, it satisfies
\begin{align}\label{S3.93}
\|wf^\v\|_{L^\infty_{x,v}} +|wf^\v|_{L^\infty(\gamma)} \leq C_{\v,d} \|\nu^{-1}wg\|_{L^\infty_{x,v}},
\end{align}
where the positive constant $C_{\v}>0$ depends only on $\v$. Moreover, if
$g$ is continuous in $\Omega_d\times \mathbb{R}^3$, then $f^\v$ is continuous away from the  grazing set $\gamma_0$.
\end{lemma}

\noindent{\bf Proof.} Let $f^n$ be the solution of \eqref{S3.54} constructed in Lemma \ref{lemS3.4} for $n\geq n_0$ with $n_0$ large enough.   Multiplying \eqref{S3.54} by $f^n$, one obtains that
\begin{align}\label{S3.94}
\v\|f^n\|^2_{L^2_{x,v}}+c_0 \|{(\mathbf{I-P})}f^n\|^2_{\nu}
\leq C_{\v} \|g\|^2_{L^2_{x,v}}.
\end{align}
We apply \eqref{3.14} and use \eqref{S3.94} to obtain
\begin{align}\label{S3.101-1}
\|wf^n\|_{L^\infty_{x,v}}{+|wf^n|_{L^\infty(\g)}}&\leq C\Big\{\|\nu^{-1}wg\|_{L^\infty_{x,v}}+\|f^n\|_{L^2_{x,v}}\Big\}\leq C_{\v,d}\|\nu^{-1}w g\|_{L^\infty_{x,v}} .
\end{align}
Taking the difference $f^{n_1}-f^{n_2}$ with $n_1,n_2\geq n_0$, we know that
\begin{align}\label{S3.102}
\begin{cases}
\v (f^{n_1}-f^{n_2})+ v_3\partial_x (f^{n_1}-f^{n_2})+\FL (f^{n_1}-f^{n_2})=0,\\[2mm]
 (f^{n_1}-f^{n_2})(x,v)|_{\gamma_{-}}=(1-\frac{1}{n_1})  (f^{n_1}-f^{n_2})(x,R_xv)+(\frac{1}{n_2}-\frac{1}{n_1}) f^{n_2}(x,R_xv).
\end{cases}
\end{align}
Multiplying \eqref{S3.102} by $ f^{n_1}-f^{n_2}$, and integrating it over $\Omega_d\times\mathbb{R}^3$, we can obtain
\begin{align}\label{S3.103}
&\v\|(f^{n_1}-f^{n_2})\|^2_{L^2_{x,v}} +c_0 \|{(\mathbf{I-P})}(f^{n_1}-f^{n_2})\|^2_{\nu}\nonumber\\
&\leq C (\frac1{n_1}+\frac1{n_2})\int_{\gamma_-}  |v_3(|f^{n_1}|+|f^{n_2}|) f^{n_2}|dv\nonumber\\
&\leq C (\frac1{n_1}+\frac1{n_2}) (|f^{n_1}|_{L^2(\gamma_+)}+|f^{n_2}|_{L^2(\gamma_+)})\nonumber\\
&\leq C_{\v,d}(\frac1{n_1}+\frac1{n_2}) \|\nu^{-1}wg\|_{L^\infty_{x,v}}\rightarrow0
\end{align}
as $n_1$, $n_2 \rightarrow \infty$,
where we have used the uniform estimate \eqref{S3.101-1} in the last inequality.  Applying \eqref{3.14} to $  f^{n_1}-f^{n_2}$  and using \eqref{S3.103}, then one has
\begin{align*}
&\|w(f^{n_1}-f^{n_2})\|_{L^\infty_{x,v}}+|w(f^{n_1}-f^{n_2})|_{L^\infty(\g)}\nonumber\\
&\leq C \Big|(\frac{1}{n_2}-\frac{1}{n_1})  wf^{n_2}\Big|_{L^\infty(\gamma_+)}+C\|f^{n_1}-f^{n_2}\|_{L^2_{x,v}}
\nonumber\\
&\leq C_{\v,d}(\frac1{n_1}+\frac1{n_2}) \|\nu^{-1}wg\|_{L^\infty_{x,v}}\rightarrow0,
\end{align*}
as $n_1,\ n_2 \rightarrow \infty$, which yields that $wf^n$ is a Cauchy sequence in $L^\infty$. We denote $f^\v=\lim_{n\rightarrow\infty} f^n$, then it is direct to check that $f^\v$ is a solution to \eqref{S3.3}, and \eqref{S3.93} holds. The continuity of $f^\v$ is easy to obtain since the convergence of  sequences is always in $L^\infty_{x,v}$ and $f^n$ is continuous away from the grazing set.  Therefore we have completed the proof of Lemma \ref{lemS3.6}. $\hfill\Box$

\

From \eqref{3.2}, \eqref{1.12}, it is direct to check that  the source term $g$ in \eqref{S3.3} satisfies
\begin{align}\label{3.57-1}
\int_{\R^3} (1,v_1-\fu_1,v_2-\fu_2,|v-\fu|^2-3) \sqrt{\mu} g(x,v) dv =(0,0,0,0).
\end{align}
Multiplying \eqref{S3.3} by $\sqrt{\mu}$ and integrating over $[0,d]\times \R^3$ to obtain
\begin{equation}\label{3.57}
\int_0^d\int_{\R^3} \sqrt{\mu} f^\v(x,v) dvdx=-\int_{\R^3} v_3\sqrt{\mu} f^\v(d,v) dv+\int_{\R^3} v_3\sqrt{\mu} f^\v(0,v) dv=0,
\end{equation}
where we have used the property of specular boundary condition in the second equality.

\smallskip

Similarly, multiplying \eqref{S3.3} by $(v_1-\fu_1,v_2-\fu_2,|v-\fu|^2-3) \sqrt{\mu}$, respectively, one gets
\begin{align}\label{3.59}
\int_0^d\int_{\R^3} (v_1-\fu_1,v_2-\fu_2,|v-\fu|^2-3) \sqrt{\mu} f^\v(x,v) dvdx=0.
\end{align}

We denote
\begin{equation*}
\mathbf{P}f^\v(x,v)=\{a^\v(x)+b^\v\cdot (v-\fu)+c^\v(x) (|v-\fu|^2-3)\}\sqrt{\mu},
\end{equation*}
then it follows from \eqref{3.57} and \eqref{3.59} that
\begin{equation}\label{3.60}
\int_0^d a^\v(x) dx=\int_0^d b^\v_1(x) dx=\int_0^d b^\v_2(x) dx=\int_0^d c^\v(x) dx=0.
\end{equation}

\begin{lemma}\label{lemS3.5}
Let $d\geq 1$. Assume \eqref{3.57-1} and  let $f^\v$ be the solution of \eqref{S3.3} constructed in Lemma \ref{lemS3.4},  then it holds that
\begin{align}\label{S3.92}
	\|\mathbf{P}f^\v\|^2_{L^2_{x,v}}\leq C d^6 \Big\{ \|{(\mathbf{I-P})}f^\v\|^2_{\nu}+\|g\|_{L^2_{x,v}}^2\Big\}.
\end{align}
\end{lemma}

\begin{remark}
By choosing suitable test function, Yin-Zhao \cite{Yin-Zhao} and Guo-Hwang-Jang-Ouyang \cite{Guo-2020} proved similar result for the time-dependent Boltzmann equation and Landau equation in three dimensional case, respectively.  The Lemma \ref{lemS3.5} can be regarded as a one dimensional version of \cite{Yin-Zhao,Guo-2020}.  We point out that   the condition \eqref{3.57-1} or \eqref{3.60} is necessary when we construct the test function in one dimensional cases.
\end{remark}

\noindent{\bf Proof.} The weak formulation of \eqref{S3.3} is
\begin{align}\label{3.62}
&\v\int_0^d\int_{\R^3} f^\v(x,v) \psi(x,v) dvdx-\int_0^d\int_{\R^3}  v_3 f^\v(x,v)\partial_x \psi(x,v)  dvdx\nonumber\\
&=-\int_{\R^3} v_3f^\v(d,v) \psi(d,v) dv+\int_{\R^3} v_3f^\v(0,v) \psi(0,v) dv\nonumber\\
&\qquad-\int_0^d\int_{\R^3} \FL f^\v(x,v) \psi(x,v) dvdx+\int_0^d\int_{\R^3} g(x,v) \psi(x,v) dvdx.
\end{align}
Similar as in \cite{Yin-Zhao} we choose some special test function $\psi$ to calculate the macroscopic part of $f^\v$.

\vspace{1mm}

\noindent{\it Step 1. Estimate on $c^\v$.} Define
\begin{equation*}
\zeta_c(x)=-\int_0^x c^\v(z)dz.
\end{equation*}
Noting \eqref{3.60}, it is easy to check that
\begin{equation}\label{3.64}
\zeta_c(x)|_{x=0,d}=0,\quad \mbox{and}\quad \|\zeta_c\|_{L^2_x}\leq d \|c^\v\|_{L^2_x}.
\end{equation}
We define the test function $\psi$ in  \eqref{3.62} to be
\begin{equation*}
\psi=\psi_c(x,v)=v_3 (|v-\fu|^2 -5) \sqrt{\mu} \zeta_c(x).
\end{equation*}
Then the second term on LHS of \eqref{3.62} is estimated as
\begin{align}\label{3.66}
&-\int_0^d\int_{\R^3}  v_3 f(x,v)\partial_x \psi_c(x,v)  dvdx\nonumber\\
&=\int_0^d\int_{\R^3}   [a^\v+b^\v\cdot (v-\fu)+c^\v(x)(|v-\fu|^2-3)] v_3^2 (|v-\fu|^2 -5) \mu(v) c^\v(x) dvdx\nonumber\\
&\qquad-\int_0^d\int_{\R^3}    (\mathbf{I-P})f^\v(x,v) v_3^2 (|v-\fu|^2 -5) \sqrt{\mu} c^\v(x) dvdx\nonumber\\
&\geq 10 \|c^\v\|_{L^2_x}^2-C\|(\mathbf{I-P})f^\v\|_{\nu}\|c^\v\|_{L^2_x}\geq 5 \|c^\v\|_{L^2_x}^2-C\|(\mathbf{I-P})f^\v\|_{\nu}^2,
\end{align}
where we have used
\begin{equation}\label{3.67}
\int_{\R^3} (|v-\fu|^2-3) v_3^2 (|v-\fu|^2 -5) \mu(v) dv=10, \quad
\int_{\R^3} v_3^2 (|v-\fu|^2 -5) \mu(v) dv=0.
\end{equation}

By using \eqref{3.67}, the first term on LHS of \eqref{3.62} is bounded as
\begin{align*}
\v\left|\int_0^d\int_{\R^3} f^\v(x,v) \psi_c(x,v) dvdx\right|
&\leq C\v\|(\mathbf{I-P})f^\v\|_{\nu} \|\zeta_c\|_{L^2_x} 
\leq C\v d \|(\mathbf{I-P})f^\v\|_{\nu} \|c^\v\|_{L^2_x}.
\end{align*}
Noting \eqref{3.64}, it is direct to have that
\begin{equation*}
-\int_{\R^3} v_3f^\v(d,v) \psi_c(d,v) dv+\int_{\R^3} v_3f^\v(0,v) \psi_c(0,v) dv=0.
\end{equation*}
Hence the RHS of \eqref{3.62} is bounded by
\begin{align}\label{3.70}
\mbox{RHS of } \eqref{3.62} \leq C d \Big( \|(\mathbf{I-P})f^\v\|_{\nu} +\|g\|_{L^2_{x,v}}\Big)\|c^\v\|_{L^2_x}.
\end{align}
Combining \eqref{3.66}-\eqref{3.70}, one obtains
\begin{align}\label{3.71}
\|c^\v\|_{L^2_x}^2\leq C d^2 \Big( \|(\mathbf{I-P})f^\v\|_{\nu}^2 +\|g\|_{L^2_{x,v}}^2\Big).
\end{align}

\vspace{1mm}

\noindent{\it Step 2. Estimate on $b^\v$.} We define
\begin{align*}
\zeta_{b,i}(x)=-\int_0^x b^\v_{i}(z)dz,\quad i=1,2,3.
\end{align*}
By using  \eqref{3.60}, it holds that
\begin{equation}\label{3.72}
\begin{split}
&\zeta_{b,i}(x)|_{x=0,d}=0, \quad i=1,2,\quad \mbox{and}\quad \zeta_{b,3}(0)=0,\\
&\|\zeta_{b,i}\|_{L^2_x}\leq d \|b^\v_i\|_{L^2_x}, \quad i=1,2,3.
\end{split}
\end{equation}
Now we take the test function $\psi$ in \eqref{3.62} to be
\begin{align*}
\psi=\psi_b(x,v)=\sum_{i=1}^3 \zeta_{b,i}(x) (v_i-\fu_i) v_3 \sqrt{\mu}-\frac12\zeta_{b,3}(x) (|v-\fu|^2-1) \sqrt{\mu}.
\end{align*}
Then the second term on LHS of \eqref{3.62} is controlled as
\begin{align}\label{3.75}
&-\int_0^d\int_{\R^3}  v_3 f^\v(x,v)\partial_x \psi_b(x,v)  dvdx\nonumber\\
&=\int_0^d\int_{\R^3}   [a^\v+b^\v\cdot (v-\fu)+c^\v(x)(|v-\fu|^2-3)] \nonumber\\
&\qquad\quad\times \left[\sum_{i=1}^3 b^\v_i(x) (v_i-\fu_i) v_3^2 -\frac12b^\v_3(x) (|v-\fu|^2-1)v_3 \right]\mu(v) dvdx\nonumber\\
&\qquad-\int_0^d\int_{\R^3}  v_3 (\mathbf{I-P})f^\v(x,v)\partial_x \psi_b(x,v)  dvdx\nonumber\\
&=\int_0^d\int_{\R^3}    \sum_{i=1}^3 |b^\v_i(x)(v_i-\fu_i)|^2 v_3^2 \mu(v)-\frac12 |b^\v_3(x)|^2 (|v-\fu|^2-1)v_3^2 \mu(v) \, dvdx\nonumber\\
&\qquad-\int_0^d\int_{\R^3}  v_3 (\mathbf{I-P})f^\v(x,v)\partial_x \psi_b(x,v)  dvdx\nonumber\\
&\geq \|(b_1^\v,b_2^\v,b_3^\v)\|_{L^2_x}^2-C\|(\mathbf{I-P})f^\v\|_{\nu} \|(b_1^\v,b_2^\v,b_3^\v)\|_{L^2_x}\nonumber\\
&\geq \frac34 \|(b_1^\v,b_2^\v,b_3^\v)\|_{L^2_x}^2-C\|(\mathbf{I-P})f^\v\|_{\nu}^2,
\end{align}
where we have used
\begin{equation*}
\begin{split}
\int_{\R^3} (v_i-\fu_i)^2 v_3^2\mu(v)dv=1,\,\, i=1,2, \quad \mbox{and}\quad
\int_{\R^3} v_3^4\mu(v)-\frac12v_3^2(|v-\fu|^2-1)\mu(v) dv=1.
\end{split}
\end{equation*}
For the first term on LHS of \eqref{3.62} is bounded as
\begin{align*}
&\v\left|\int_0^d\int_{\R^3} f^\v(x,v) \psi_b(x,v) dvdx\right| \nonumber\\
&\leq C\v d \|c^\v\|_{L^2_x}\|(b_1^\v,b_2^\v,b_3^\v)\|_{L^2_x}+Cd\v\|(\mathbf{I-P})f^\v\|_{\nu} \|(b_1^\v,b_2^\v,b_3^\v)\|_{L^2_x} \nonumber\\
&\leq \frac14 \|(b_1^\v,b_2^\v,b_3^\v)\|_{L^2_x}^2+C\v^2 d^2 \Big(\|(\mathbf{I-P})f^\v\|_{\nu}^2+ \|c^\v\|_{L^2_x}^2\Big).
\end{align*}
For the boundary terms on RHS of \eqref{3.62}, it follows from \eqref{3.72} that
\begin{align*}
&-\int_{\R^3} v_3f^\v(d,v) \psi_b(d,v) dv+\int_{\R^3} v_3f^\v(0,v) \psi_b(0,v) dv\nonumber\\
&=-\zeta_{b,3}(d)\cdot\int_{\R^3} v_3f^\v(d,v)\cdot  [v_3^2-\frac12(|v-\fu|^2-1)] \sqrt{\mu}dv=0,
\end{align*}
due to  the fact that $f(d,v)$ is an even function with respect to $v_3$. Thus the terms on RHS of \eqref{3.62}, we have
\begin{equation}\label{3.79}
\mbox{RHS of } \eqref{3.62} \leq Cd \|(b_1^\v,b_2^\v,b_3^\v)\|_{L^2_x} \Big\{\|(\mathbf{I-P})f^\v\|_{\nu}+ \|g\|_{L^2_{x,v}}\Big\}.
\end{equation}
Then combining \eqref{3.75}-\eqref{3.79}, one obtains that
\begin{align}\label{3.80}
\|(b_1^\v,b_2^\v,b_3^\v)\|_{L^2_x}^2\leq Cd^4\Big\{\|(\mathbf{I-P})f^\v\|_{\nu}^2+ \|g\|_{L^2_{x,v}}^2\Big\}.
\end{align}

\vspace{1mm}

\noindent{\it Step 2. Estimate on $a^\v$.} Define
\begin{equation*}
\zeta_a(x)=-\int_0^x a^\v(z)dz.
\end{equation*}
Noting \eqref{3.60}, it is easy to check that
\begin{equation}\label{3.81}
\zeta_a(x)|_{x=0,d}=0,\quad \mbox{and}\quad \|\zeta_a\|_{L^2}\leq d \|a^\v\|_{L^2_x}.
\end{equation}
We define the test function $\psi$ in  \eqref{3.62} to be
\begin{equation*}
\psi=\psi_a(x,v)=(|v-\fu|^2 -10) \sqrt{\mu}\, v_3\, \zeta_a(x).
\end{equation*}
Then the second term on LHS of \eqref{3.62} is estimated as
\begin{align}\label{3.84}
&-\int_0^d\int_{\R^3}  v_3 f^\v(x,v)\partial_x \psi_a(x,v)  dvdx\nonumber\\
&=\int_0^d\int_{\R^3}   [a^\v+b^\v\cdot (v-\fu)+c^\v(x)(|v-\fu|^2-3)]\cdot a^\v(x)(|v-\fu|^2 -10) v_3 ^2\mu(v) \nonumber\\
&\qquad-\int_0^d\int_{\R^3}  v_3 (\mathbf{I-P})f^\v(x,v)\partial_x \psi_a(x,v)  dvdx\nonumber\\
&\geq 5\|a^\v\|_{L^2_x}^2-C\|(\mathbf{I-P})f^\v\|_{\nu}\|a^\v\|_{L^2} \geq 4\|a^\v\|_{L^2_x}^2-C\|(\mathbf{I-P})f^\v\|_{\nu}^2,
\end{align}
where we have used
\begin{align}\nonumber
\int_{\R^3} (|v-\fu|^2-3) v_3^2 \cdot (|v-\fu|^2 -10) \mu(v) dv=0.
\end{align}
A direct calculation shows that
\begin{align*}
&\v\left|\int_0^d\int_{\R^3} f^\v(x,v) \psi_a(x,v) dvdx\right| \nonumber\\
&\leq \left|\v\int_0^d\int_{\R^3} a^\v(x) \zeta_{a}^\v(x) v_3^2(|v-\fu|^2 -10) \mu(v)dvdx \right|+C\v d\|(\mathbf{I-P})f^\v\|_{\nu} \|(a^\v)\|_{L^2_x} \nonumber\\
&\leq \|a^\v\|_{L^2_x}^2+Cd^2\v^2 (\|(\mathbf{I-P})f^\v\|_{\nu}+\|b^\v_3\|_{L^2_x}^2).
\end{align*}
By using \eqref{3.81}, one has
\begin{equation*}
-\int_{\R^3} v_3f^\v(d,v) \psi_a(d,v) dv+\int_{\R^3} v_3f^\v(0,v) \psi_a(0,v) dv=0.
\end{equation*}
Hence, for the RHS of \eqref{3.62}, it holds that
\begin{align}\label{3.87}
\mbox{RHS of } \eqref{3.62} \leq Cd \|a^\v\|_{L^2_x} \Big\{\|(\mathbf{I-P})f^\v\|_{\nu}+ \|g\|_{L^2_{x,v}}\Big\}.
\end{align}
Combining \eqref{3.84}-\eqref{3.87} and using \eqref{3.80}, we obtain
\begin{equation}\label{3.88}
\|a^\v\|_{L^2_x}^2\leq Cd^6 \{\|(\mathbf{I-P})f^\v\|_{\nu}^2+ \|g\|_{L^2_x}^2\}.
\end{equation}
Therefore, \eqref{S3.92} follows directly from \eqref{3.88}, \eqref{3.80} and \eqref{3.71}. $\hfill\Box$

\vspace{2mm}

\begin{lemma}\label{lemS3.7}
Let $d\geq 1$, $\beta\geq 3$.  Assume \eqref{3.57-1} and  $\|\nu^{-1}wg\|_{L^\infty_{x,v}}<\infty$.   Then there exists a  unique solution $f=f(x,v)$  to the linearized steady Boltzmann equation
\begin{align}\label{S3.106}
\begin{cases}
v_3 \partial_xf+\FL f=g,\quad (x,v)\in \Omega_d\times \R^3, \\
f(x,v)|_{\gamma_{-}}=f(x,R_xv),
\end{cases}
\end{align}
with
\begin{equation}\label{S3.107}
\|wf\|_{L^\infty_{x,v}} +|wf|_{L^\infty(\gamma)} \leq C_d \|\nu^{-1}wg\|_{L^\infty_{x,v}}.
\end{equation}
Moreover, if $g$ is continuous in $\Omega_d\times\mathbb{R}$,
then $f$ is continuous away from the grazing set $\gamma_0$.
\end{lemma}

\noindent{\bf Proof.} Let $f^\v$ be the solution of \eqref{S3.3} constructed in Lemma \ref{lemS3.6} for $\v>0$.   Multiplying the first equation of $\eqref{S3.3}_1$ by $f^\v$ and integrating the resultant equation over $\Omega_d\times\mathbb{R}^3$, we have
\begin{align}\label{S3.109}
\v\|f^\v\|^2_{L^2_{x,v}}+c_0\|(\mathbf{I-P})f^\v\|^2_{\nu}\leq \delta \|f^\v\|^2_{L^2_{x,v}}+C_\delta \|g\|^2_{L^2_{x,v}}.
\end{align}
which, together with  Lemma \ref{lemS3.5}, yields that
\begin{equation*}
\|f^\v\|^2_{L^2_{x,v}} \leq C d^6 \delta \|f^\v\|^2_{L^2_{x,v}}+C_{d, \delta} \|g\|^2_{L^2_{x,v}}.
\end{equation*}
Taking $\delta$ small enough, we obtain
\begin{equation}\label{3.94}
\|f^\v\|^2_{L^2_{x,v}} \leq C_{d} \|g\|^2_{L^2_{x,v}}.
\end{equation}

Applying \eqref{3.14} to $f^\v$ and using \eqref{3.94}, then we obtain
\begin{equation}\label{S3.119}
\|wf^\v\|_{L^\infty_{x,v}}{+|wf^\vep|_{L^\infty(\g)}}\leq C_d\|\nu^{-1}wg\|_{L^\infty_{x,v}}.
\end{equation}
Next we consider the convergence of $f^\v$ as $\v\rightarrow0+$. For any $\v_1,\v_2>0$, we consider the difference $f^{\v_2}-f^{\v_1}$ satisfying
\begin{equation}\label{S3.120}
\begin{cases}
v_3\partial_x (f^{\v_2}-f^{\v_1})+ \FL(f^{\v_2}-f^{\v_1})=-\v_2 f^{\v_2}+\v_1 f^{\v_1},\\[2mm]
(f^{\v_2}-f^{\v_1})|_{\gamma_-}=(f^{\v_2}-f^{\v_1})(x,R_xv).
\end{cases}
\end{equation}
Multiplying \eqref{S3.120} by $f^{\v_2}-f^{\v_1}$, integrating the resultant equation and by similar arguments as in \eqref{S3.109}-\eqref{3.94}, one gets
\begin{align}\label{S3.121}
\|f^{\v_2}-f^{\v_1}\|^2_{L^2_{x,v}}
&\leq C_d \|\v_2 f^{\v_2}-\v_1 f^{\v_1}\|^2_{L^2_{x,v}}
\leq C_d(\v_1^2+\v_2^2) \|g\|^2_{L^2_{x,v}}\nonumber\\
&\leq C_d (\v_1^2+\v_2^2)\cdot  \|\nu^{-1}wg\|_{L^\infty_{x,v}}^2\rightarrow0,
\end{align}
as $\v_1$, $\v_2\rightarrow 0+$. Finally, applying \eqref{3.14} to $f^{\v_2}-f^{\v_1}$ and using \eqref{S3.121}, then we obtain
\begin{align}\label{S3.122}
&\|w(f^{\v_2}-f^{\v_1})\|_{L^\infty_{x,v}}{+| w(f^{\v_2}-f^{\v_1})|_{L^\infty(\g)}}\nonumber\\
&\leq C\Big\{ \|\nu^{-1}w (\v_2 f^{\v_2}-\v_1 f^{\v_1})\|_{L^\infty_{x,v}} +\|f^{\v_2}-f^{\v_1}\|_{L^2_{x,v}} \Big\}\nonumber\\
&\leq C_d(\v_1+\v_2)\|\nu^{-1}wg\|_{L^\infty_{x,v}}\rightarrow 0,
\end{align}
as $\v_1$, $\v_2\rightarrow 0+$,
 With \eqref{S3.122}, we know that there exists a function $f$ so that $\|w(f^{\v}-f)\|_{L^\infty_{x,v}}\rightarrow0$ as $\v\rightarrow 0+$. And it is  direct to see that $f$ solves \eqref{S3.106}. Also,  \eqref{S3.107} follows immediately from \eqref{S3.119}. The continuity of $f$
directly follows from the $L^\infty_{x,v}$-convergence and the continuity of $f^\v$. Therefore the proof of Lemma \ref{lemS3.7} is complete. $\hfill\Box$

\

To obtain the solution for half-space problem, we need some uniform estimate independent of $d$, then we can take the limit $d\rightarrow\infty$. Let $f$ be the solution of \eqref{S3.106}, we denote
\begin{align*}
\mathbf{P} f(x,v)=\big[a(x)+b(x)\cdot (v-\fu)+c(x) (|v-\fu|^2-3)\big] \sqrt{\mu}.
\end{align*}
Multiplying \eqref{S3.106} by $\sqrt{\mu}$ and using \eqref{3.57-1}, we have
\begin{align}\label{3.100}
0=\frac{d}{dx} \int_{\R^3} v_3 \sqrt{\mu} f(x,v)dv = \frac{d}{dx} b_3(x)\equiv0.
\end{align}
Since $f$ satisfies the specular boundary, it holds that $b_3(x)|_{x=0}=b_3(x)|_{x=d}=0$, which, together with \eqref{3.100} yields
\begin{equation}\label{3.101}
b_3(x)=0,\quad \mbox{for}\  x\in[0,d].
\end{equation}

Let $(\phi_0,\phi_1,\phi_2,\phi_3)$ be some constants chosen later, we define
\begin{align*}
\bar{f}(x,v)&:=f(x,v)+[\phi_0+\phi_1 (v_1-\fu_1)+\phi_2 (v_2-\fu_2)+\phi_3 (|v-\fu|^2-3)]\sqrt{\mu}\nonumber\\
&=[\bar{a}(x)+\bar{b}_1(x)\cdot (v_1-\fu_1)+\bar{b}_2(x)\cdot (v_2-\fu_2)+\bar{c}(x) (|v-\fu|^2-3)]\sqrt{\mu}\nonumber\\
&\qquad +(\mathbf{I-P}) \bar{f},
\end{align*}
where
\begin{align*}
\begin{cases}
\bar{a}(x)=a(x)+\phi_0,\\
\bar{b}_i(x)=b_i(x)+\phi_i, \quad i=1,2,\\
\bar{c}(x)=c(x)+\phi_3.
\end{cases}
\end{align*}
From \eqref{3.101}, it is easy to check that
\begin{equation}\label{3.104-1}
\bar{b}_3(x)=0\quad \mbox{and}\quad(\mathbf{I-P}) \bar{f}(x,v)=(\mathbf{I-P}) f(x,v)\quad \forall x\in[0,d].
\end{equation}
In fact, it is direct to check that $\bar{f}$ still satisfies \eqref{S3.106}, i.e.,
\begin{align}\label{3.104}
\begin{cases}
v_3 \partial_x\bar{f}+\FL \bar{f}=g,\quad  (x,v)\in\R^3,\\
\bar{f}(x,v)|_{\gamma_{-}}=\bar{f}(x,R_xv),
\end{cases}
\end{align}
with
\begin{equation*}
\|w\bar{f}\|_{L^\infty_{x,v}} +|w\bar{f}|_{L^\infty(\gamma)} \leq C_d \|wg\|_{L^\infty_{x,v}}+C_d |(\phi_0,\phi_1,\phi_2,\phi_3)|.
\end{equation*}
Multiplying \eqref{3.104} by $(v_1-\fu_1,v_2-\fu_2,|v-\fu|^2-5) \sqrt{\mu}$ and using \eqref{3.57-1}, we get
\begin{align}\label{3.102}
\begin{split}
\int_{\R^3} v_3 (v_i-\fu_i) \sqrt{\mu} \bar{f}(x,v)dv&=0, \quad \forall\, x\in[0,d],\quad i=1,2,\\
\int_{\R^3} v_3 (|v-\fu|^2-5) \sqrt{\mu} \bar{f}(x,v)dv&=0,\quad \forall \, x\in[0,d].
\end{split}
\end{align}

It follows from \eqref{3.104-1} and \eqref{3.102} that
\begin{align}\label{3.108}
&\int_{\R^3} v_3|\mathbf{P}\bar{f}(x,v)|^2 dv\nonumber\\
&=\int_{\R^3} v_3 |[\bar{a}+\bar{b}_1\cdot (v_1-\fu_1)+\bar{b}_2\cdot (v_2-\fu_2)+\bar{c} (|v-\fu|^2-3)]|^2\mu(v)dv\equiv0,
\end{align}
and
\begin{align}\label{3.109}
&\int_{\R^3} v_3\mathbf{P}\bar{f}(x,v)\cdot (\mathbf{I-P})\bar{f}(x,v) dv\nonumber\\
&=\int_{\R^3} v_3 [\bar{a}+\bar{b}_1\cdot (v_1-\fu_1)+\bar{b}_2\cdot (v_2-\fu_2)+\bar{c} (|v-\fu|^2-3)]\sqrt{\mu(v)} \cdot (\mathbf{I-P})\bar{f}  dv \equiv0.
\end{align}
By utilizing \eqref{3.108} and \eqref{3.109}, it holds that
\begin{align}\label{3.110}
\int_{\R^3} v_3|\bar{f}(x,v)|^2 dv=\int_{\R^3} v_3|(\mathbf{I-P})\bar{f}(x,v)|^2 dv,\quad \forall x\in[0,d].
\end{align}

Multiplying \eqref{3.104} by $\bar{f}$ and using \eqref{3.110}, \eqref{3.57-1}, we obtain
\begin{align}\label{3.111}
\frac{d}{dx} \int_{\R^3} v_3 |(\mathbf{I-P})\bar{f}|^2dv+c_0\|(\mathbf{I-P})\bar{f}\|^2_{\nu}\leq C \|g\|_{L^2_v}^2,
\end{align}
where we have used the fact
\begin{align}\nonumber
\int_{\R^3} g \bar{f}dv=\int_{\R^3} g (\mathbf{I-P})\bar{f} dv\leq \frac12c_0\|(\mathbf{I-P})\bar{f}\|_{\nu}^2+C\|g\|_{L^2_v}^2.
\end{align}

Let $0\leq \sigma<\sigma_1\leq \sigma_0$. Multiplying \eqref{3.111} by $e^{2\s_1 x}$, one obtains that
\begin{align*}
\frac{d}{dx}\left\{e^{2\s_1 x}\int_{\R^3} v_3 |(\mathbf{I-P})\bar{f}|^2dv \right\} +(c_0-C\s_1)e^{2\s_1 x}  \|(\mathbf{I-P})\bar{f}\|^2_{\nu}
\leq Ce^{2\s_1 x} \|g\|_{L^2_v}^2.
\end{align*}
Taking $\s_0>0$ small such that $c_0-C\s_1\geq\frac12 c_0$, then we have
\begin{align}\label{3.113}
\frac{d}{dx}\left\{e^{2\s_1 x}\int_{\R^3} v_3 |(\mathbf{I-P})\bar{f}|^2dv \right\} +\frac12c_0e^{2\s_1 x}  \|(\mathbf{I-P})\bar{f}\|^2_{\nu}
\leq Ce^{2\s_1 x} \|g\|_{L^2_v}^2.
\end{align}
Integrating \eqref{3.113} over $[0,d]$ and noting $\eqref{3.104-1}_2$ one has
\begin{equation}\label{3.114}
\int_0^d e^{2\s_1 x}  \|(\mathbf{I-P})f\|^2_{\nu} dx\equiv \int_0^d e^{2\s_1 x}  \|(\mathbf{I-P})\bar{f}\|^2_{\nu} dx
\leq C\|e^{\s_1 x}g\|_{L^2_{x,v}}^2,
\end{equation}
where we have used the fact $\dis\int_{\R^3} v_3 |(\mathbf{I-P})\bar{f}(0,v)|^2dv=\int_{\R^3} v_3 |(\mathbf{I-P})\bar{f}(d,v)|^2dv$ due to the specular boundary condition.

\begin{remark}
We point out that the estimations \eqref{3.102}-\eqref{3.114} are independent of the choice of $(\phi_0,\phi_1,\phi_2,\phi_3)$. However, to obtain uniform estimate independent of $d$ for macroscopic part, we need  to choose  $(\phi_0,\phi_1,\phi_2,\phi_3)$ suitably.
\end{remark}

We denote
\begin{align}\label{3.115}
\begin{split}
\mathcal{A}_{ij}(v)&=\{(v_i-\fu_i)(v_j-\fu_j)-\frac{\d_{ij}}3 |v-\fu|^2\}\sqrt{\mu},\\
\mathcal{B}_{i}(v)&= (v_i-\fu_i) (|v-\fu|^2-5)\sqrt{\mu},
\end{split}
\quad i,j=1,2,3.
\end{align}
It is obviously that $\mathcal{A}_{ij}, \mathcal{B}_{i} \in \mathbb{N}^{\perp}$, and
\begin{align}
\begin{split}
\kappa_1:=&\int_{\R^3} \mathcal{A}_{31} \FL^{-1} \mathcal{A}_{31} dv\equiv\int_{\R^3} \mathcal{A}_{ij} \FL^{-1} \mathcal{A}_{ij} dv>0,\ \mbox{for}\  i\neq j,\\
\kappa_2:=&\int_{\R^3} \mathcal{B}_{3} \FL^{-1} \mathcal{B}_{3} dv
\equiv \int_{\R^3}\mathcal{B}_{i} \FL^{-1} \mathcal{B}_{i} dv>0,\  \mbox{for}\  i=2,3.
\end{split}
\end{align}

\begin{lemma}\label{lem3.10}
There exist constants $(\phi_0,\phi_1,\phi_2,\phi_3)$  such that
\begin{align}\label{3.118}
\begin{split}
&\int_{\R^3} v_3\bar{f}(d,v)\cdot v_3 \sqrt{\mu} dv=0,\\
&\int_{\R^3} v_3\bar{f}(d,v)\cdot \FL^{-1} (\mathcal{A}_{3i}) dv=0,\  i=1,2,\\
&\int_{\R^3} v_3\bar{f}(d,v)\cdot \FL^{-1} (\mathcal{B}_{3}) dv=0.
\end{split}
\end{align}
\end{lemma}

\noindent{\bf Proof.} A direct calculation shows that \eqref{3.118} is equivalent to the following
\begin{align}
\int_{\R^3} v_3\bar{f}(x,v)\cdot v_3 \sqrt{\mu} dv&=\bar{a}(x) +2\bar{c}(x) +\int_{\R^3} \mathcal{A}_{33}(v)\cdot (\mathbf{I-P})\bar{f}(x,v) dv\nonumber\\
&=\phi_0+2\phi_3 +a(x)+2c(x)+\int_{\R^3} \mathcal{A}_{33}(v)\cdot (\mathbf{I-P})f(x,v) dv,\label{3.119}\\
\int_{\R^3} v_3\bar{f}(x,v)\cdot \FL^{-1} (\mathcal{A}_{31}) dv
&= \kappa_1 \bar{b}_1(x)+\int_{\R^3} v_3 (\mathbf{I-P})\bar{f}(x,v)\cdot \FL^{-1} (\mathcal{A}_{31})  dv\nonumber\\
&=\kappa_1 \phi_1+\kappa_1 b_1(x) +\int_{\R^3} v_3 (\mathbf{I-P})f(x,v)\cdot \FL^{-1} (\mathcal{A}_{31})  dv\label{3.120} \\
\int_{\R^3} v_3\bar{f}(x,v)\cdot \FL^{-1} (\mathcal{A}_{32}) dv
&=\kappa_1 \bar{b}_2(x)+\int_{\R^3} v_3 (\mathbf{I-P})\bar{f}(x,v)\cdot \FL^{-1} (\mathcal{A}_{32})  dv\nonumber\\
&=\kappa_1 \phi_2+\kappa_1 b_2(x) +\int_{\R^3} v_3 (\mathbf{I-P})f(x,v)\cdot \FL^{-1} (\mathcal{A}_{32})  dv,\label{3.121}\\
\int_{\R^3} v_3\bar{f}(x,v)\cdot \FL^{-1} (\mathcal{B}_{3}) dv
&=\kappa_2 \bar{c}(x)+\int_{\R^3} v_3 (\mathbf{I-P})\bar{f}(x,v)\cdot \FL^{-1} (\mathcal{B}_{3})  dv\nonumber\\
&=\kappa_2 \phi_3+\kappa_2 c(x)+\int_{\R^3} v_3 (\mathbf{I-P})f(x,v)\cdot \FL^{-1} (\mathcal{B}_{3})  dv,\label{3.122}
\end{align}
where we have used $\eqref{3.104}_2$.

Using  \eqref{3.119}-\eqref{3.122}, then \eqref{3.118} is equivalent  as
\begin{equation}\nonumber
\left(
\begin{array}{cccc}
1 & 0 & 0 & 2  \\
0 & \kappa_1 & 0 & 0  \\
0 & 0 & \kappa_1 & 0  \\
0 & 0 & 0 & \kappa_2  \\
\end{array}
\right)
\left(
\begin{array}{c}
\phi_0 \\
\phi_1\\
\phi_2\\
\phi_3\\
\end{array}
\right)
=-\left(
\begin{array}{c}
a(d)+2c(d)+\int_{\R^3}(\mathbf{I-P})f(d,v) \cdot \mathcal{A}_{33}(v)dv \\[1mm]
\kappa_1 b_1(d) +\int_{\R^3} v_3 (\mathbf{I-P})f(d,v)\cdot \FL^{-1} (\mathcal{A}_{31})  dv\\[1mm]
\kappa_1 b_2(d) +\int_{\R^3} v_3 (\mathbf{I-P})f(d,v)\cdot \FL^{-1} (\mathcal{A}_{32})  dv\\[1mm]
\kappa_2 c(d)+\int_{\R^3} v_3 (\mathbf{I-P})f(d,v)\cdot \FL^{-1} (\mathcal{B}_{3})  dv
\end{array}
\right).
\end{equation}
Noting the matrix is non-singular, hence  $(\phi_0,\phi_1,\phi_2,\phi_3)$  is found. Therefore the proof of Lemma \ref{lem3.10} is completed. $\hfill\Box$

\begin{lemma}\label{lem3.11}
Let $(\phi_0,\phi_1,\phi_2,\phi_3)$ be the ones determined  in Lemma \ref{lem3.10}, then it holds that
\begin{equation}\label{3.123}
\|e^{\s x}\bar{f}\|_{L^2_{x,v}}\leq \frac{C}{\s_1-\s}\|e^{\s_1 x} g\|_{L^2_{x,v}}
\end{equation}
with $0<\s<\s_1\leq \s_0$, and the constant $C$ is independent of $d$.
\end{lemma}

\noindent{\bf Proof.}  Multiplying \eqref{3.104} by  $\FL^{-1} (\mathcal{A}_{31}), \FL^{-1} (\mathcal{A}_{32}) $ and $\FL^{-1} (\mathcal{B}_{3}) $, respectively, then we obtain
\begin{align}
\frac{d}{dx}\left(
\begin{array}{c}
\int_{\R^3} v_3\bar{f}(x,v)\cdot \FL^{-1} (\mathcal{A}_{31}) dv\\[1mm]
\int_{\R^3} v_3\bar{f}(x,v)\cdot \FL^{-1} (\mathcal{A}_{32}) dv\\[1mm]
\int_{\R^3} v_3\bar{f}(x,v)\cdot \FL^{-1} (\mathcal{B}_{3}) dv\\
\end{array}
\right)
=
\left(
\begin{array}{c}
\int_{\R^3}  [g-\FL (\mathbf{I-P})f]\cdot \FL^{-1} (\mathcal{A}_{31}) dv\\[1mm]
\int_{\R^3} [g-\FL (\mathbf{I-P})f]\cdot \FL^{-1} (\mathcal{A}_{32}) dv\\[1mm]
\int_{\R^3} [g-\FL (\mathbf{I-P})f]\cdot \FL^{-1} (\mathcal{B}_{3}) dv\\
\end{array}
\right).\nonumber
\end{align}
Integrating above system over $[x,d]$ and using $\eqref{3.118}_{2,3,4}$ to get
\begin{align}
\left(
\begin{array}{c}
\int_{\R^3} v_3\bar{f}(x,v)\cdot \FL^{-1} (\mathcal{A}_{31}) dv\\[1mm]
\int_{\R^3} v_3\bar{f}(x,v)\cdot \FL^{-1} (\mathcal{A}_{32}) dv\\[1mm]
\int_{\R^3} v_3\bar{f}(x,v)\cdot \FL^{-1} (\mathcal{B}_{3}) dv\\
\end{array}
\right)
=
\int_x^d\left(
\begin{array}{c}
\int_{\R^3}  [\FL (\mathbf{I-P})f-g]\cdot \FL^{-1} (\mathcal{A}_{31}) dv\\[1mm]
\int_{\R^3} [\FL (\mathbf{I-P})f-g]\cdot \FL^{-1} (\mathcal{A}_{32}) dv\\[1mm]
\int_{\R^3} [\FL (\mathbf{I-P})f-g]\cdot \FL^{-1} (\mathcal{B}_{3}) dv\\
\end{array}
\right)(z)dz,\nonumber
\end{align}
which, together with \eqref{3.120}-\eqref{3.122}, yields that
\begin{align}\label{3.125}
|(\kappa_1 \bar{b}_{1},\kappa_1  \bar{b}_{2},\kappa_2 \bar{c})(x)|\leq C\|(\mathbf{I-P})f(x)\|_{\nu}+C\int_x^d\{\|(\mathbf{I-P})f(z)\|_{\nu}+\|g(z)\|_{L^2_v}\}dz.
\end{align}
Multiplying \eqref{3.125} by $e^{\s x}$ with $0<\s<\s_1\leq\s_0$ and using \eqref{3.114}, then we can have
\begin{align}\label{3.126}
\int_0^d e^{2\s x}|(\bar{b}_{1},\bar{b}_{2},\bar{c})(x)|^2d x
&\leq  C\int_0^d e^{2\s x} \|(\mathbf{I-P})f(x)\|_{\nu}^2 dx\nonumber\\
&\hspace{-6mm}+C\int_0^d e^{2\s x} \left|\int_x^d\{\|(\mathbf{I-P})f(z)\|_{\nu}+\|g(z)\|_{L^2_v}\}dz\right|^2 dx\nonumber\\
&\leq \frac{C}{\s_1-\s}\|e^{\s_1 x} g\|_{L^{2}_{x,v}}^2.
\end{align}

Finally we consider the case for $\bar{a}$. In fact, multiplying \eqref{3.104} by $v_3\sqrt{\mu}$, we have that
\begin{align}\nonumber
\frac{d}{dx} \int_{\R^3} \bar{f}(x,v)\cdot v_3^2 \sqrt{\mu} dv=\int_{\R^3}g\cdot v_3\sqrt{\mu}  dv.
\end{align}
Integrating above equation over $[x,d]$,  and using $\eqref{3.118}_1$, $\eqref{3.104-1}_2$, one obtain
\begin{align}\label{3.128}
\bar{a}(x)=-2\bar{c}(x)+\int_{\R^3}   (\mathbf{I-P})\bar{f}(x,v)\cdot v_3^2 \sqrt{\mu} dv-\int_x^d\int_{\R^3}g\cdot v_3\sqrt{\mu}  dvdx.
\end{align}
Multiplying \eqref{3.128} by $e^{\s x}$ with $0<\s<\s_1\leq \s_0$ and using \eqref{3.114}, \eqref{3.126}, then we can get
\begin{equation}\label{3.129}
\int_0^d e^{2\s x}|\bar{a}(x)|^2d x\leq  \frac{C}{\s_1-\s}\|e^{\s_1 x} g\|_{L^{2}_{x,v}}^2.
\end{equation}
Combining \eqref{3.114}, \eqref{3.126} and  \eqref{3.129}, we prove \eqref{3.123}. Therefore the proof of Lemma \ref{lem3.11} is completed. $\hfill\Box$

\

\begin{lemma}\label{lem3.12}
Let $\beta\geq 3, d\geq 1$, and $\bar{f}$ to the solution of \eqref{3.104}, it holds that
\begin{align}\label{3.129-0}
\|e^{\s x}w \bar{f}\|_{L^\infty_{x,v}} +|e^{\s x}w \bar{f}|_{L^\infty(\gamma)} \leq \frac{C}{\s_0-\s}\|e^{\s_0x}\nu^{-1} wg\|_{L^\infty_{x,v}}.
\end{align}
where the constant $C>0$ is independent of $d$.
\end{lemma}

\noindent{\bf Proof.} Let $\bar{h}:=e^{\s x} w \bar{f}$. Multiplying \eqref{3.104} by $e^{\s x}w $ to have
\begin{equation*}
v_3 \partial_x\bar{h}+\nu_{\s}(v) \bar{h}=K_w \bar{h}+e^{\s x}wg,\quad
\bar{h}(x,v)|_{\gamma_{-}}=\bar{h}(x,R_xv),
\end{equation*}
where $\nu_{\s}(v):=\nu(v) -\s v_3$. We further take $\s_0>0$ small such that $\nu_{\s}(v) \geq \frac12 \nu(v)>0$. By the same arguments as  in Lemma \ref{lemS3.3}, we can obtain
\begin{align*}
\|\bar{h}\|_{L^\infty_{x,v}} +|\bar{h}|_{L^\infty(\gamma)} & \leq C\|e^{\s x} \nu^{-1}wg\|_{L^\infty_{x,v}}+C\|e^{\s x} \bar{f}\|_{L^2_{x,v}}\nonumber\\
&\leq C\|e^{\s x} \nu^{-1}wg\|_{L^\infty_{x,v}}+\frac{C}{\s_1-\s}\|e^{\s_1 x} g\|_{L^2_{x,v}}\nonumber\\
&\leq \frac{C}{\s_0-\s}\|e^{\s_0 x} \nu^{-1}wg\|_{L^\infty_{x,v}}
\end{align*}
where we have used \eqref{3.123} and chosen $\s_1=\s+\frac{\s_0-\s}{2}$ such that  $0<\s<\s_1<\s_0$. Hence the proof of Lemma \ref{lem3.12} is completed. $\hfill\Box$

\

We shall prove Theorem \ref{thm3.1} by taking the limit $d\rightarrow\infty$. From now on, we shall denote the solution $\bar{f}(x,v)$ of \eqref{3.104} to be $\bar{f}_d(x,v)$ to emphasize the dependence of  parameter $d$. We denote
\begin{equation*}
\tilde{f}(x,v)=\bar{f}_{d_2}(x,v)-\bar{f}_{d_1}(x,v), \quad 1\leq d_1\leq d_2.
\end{equation*}
Then $\tilde{f}$ satisfies the following equation
\begin{equation}\label{3.132}
\begin{cases}
v_3\partial_x \tilde{f}+\FL \tilde{f}=0,\quad x\in[0,d_1],\ v\in\R^3, \\
\tilde{f}(0,v)|_{v_3>0}=\tilde{f}(0,Rv).
\end{cases}
\end{equation}

\subsection{Proof of Theorem \ref{thm3.1}.} As previous, we  divide the proof into two steps.\vspace{1mm}

\noindent{\it Step 1. Convergence in $L^2$-norm. } Multiplying \eqref{3.132} by $\tilde{f}$ and using \eqref{3.129-0} to obtain
\begin{align}\label{3.134}
&\int_0^{d_1}\int_{\R^3}(1+|v|)|(\mathbf{I-P})\tilde{f}(x,v)|^2 dvdx\nonumber\\
&\leq C\int_{\R^3} |v_3|\cdot |\tilde{f}(d_1,v)|^2dv \leq C \big\{\|w\bar{f}_{d_2}(d_1)\|^2_{L^\infty_{v}}+|w\bar{f}_{d_1}(d_1)|^2_{L^\infty(\gamma)}\big\}\nonumber\\
&\leq \frac{C}{(\s_0-\s)^2}\|e^{\s_0 x} \nu^{-1}wg\|^2_{L^\infty_{x,v}} e^{-2\s d_1}.
\end{align}

We still need to control the macroscopic part. We denote
\begin{align*}
\mathbf{P}\tilde{f}=[\tilde{a}(x)+\tilde{b}_{1}(x) (v_1-\fu_1)+\tilde{b}_{2}(x) (v_2-\fu_2)+ \tilde{c}(x) (|v-\fu|^2-3)]\sqrt{\mu}.
\end{align*}
Similar as in Lemma \ref{lem3.11}, multiplying \eqref{3.132} by $\FL^{-1} (\mathcal{A}_{31}), \FL^{-1} (\mathcal{A}_{32}) $ and $\FL^{-1} (\mathcal{B}_{3}) $, respectively, integrating the resultant equation over $[x,d_1]$ to have
\begin{align}
\left(
\begin{array}{c}
\int_{\R^3} v_3\tilde{f}(x,v)\cdot \FL^{-1} (\mathcal{A}_{31}) dv\\[1mm]
\int_{\R^3} v_3\tilde{f}(x,v)\cdot \FL^{-1} (\mathcal{A}_{32}) dv\\[1mm]
\int_{\R^3} v_3\tilde{f}(x,v)\cdot \FL^{-1} (\mathcal{B}_{3}) dv\\
\end{array}
\right)
&=\left(
\begin{array}{c}
\int_{\R^3} v_3\tilde{f}(d_1,v)\cdot \FL^{-1} (\mathcal{A}_{31}) dv\\[1mm]
\int_{\R^3} v_3\tilde{f}(d_1,v)\cdot \FL^{-1} (\mathcal{A}_{32}) dv\\[1mm]
\int_{\R^3} v_3\tilde{f}(d_1,v)\cdot \FL^{-1} (\mathcal{B}_{3}) dv\\
\end{array}
\right)\nonumber\\
&\quad+\int_x^{d_1}\left(
\begin{array}{c}
\int_{\R^3}  \FL (\mathbf{I-P})\tilde{f}\cdot \FL^{-1} (\mathcal{A}_{31}) dv\\[1mm]
\int_{\R^3} \FL (\mathbf{I-P})\tilde{f}\cdot \FL^{-1} (\mathcal{A}_{32}) dv\\[1mm]
\int_{\R^3} \FL (\mathbf{I-P})\tilde{f}\cdot \FL^{-1} (\mathcal{B}_{3}) dv\\
\end{array}
\right)(z)dz,\nonumber
\end{align}
which, together with \eqref{3.120}-\eqref{3.122}, yields that
\begin{align}\label{3.136}
|(\kappa_1 \tilde{b}_{1},\kappa_1  \tilde{b}_{2},\kappa_2 \tilde{c})(x)|&\leq C\Big\{\|w\bar{f}_{d_2}(d_1)\|_{L^\infty_{v}}+|w\bar{f}_{d_1}(d_1)|_{L^\infty(\g)}\Big\}+C\|(\mathbf{I-P})\tilde{f}(x)\|_{\nu}\nonumber\\
&\quad+C\int_x^{d_1}\|(\mathbf{I-P})\tilde{f}(z)\|_{\nu} dz.
\end{align}
Integrating \eqref{3.136} over $[0,d_1]$ and using \eqref{3.129-0}, \eqref{3.134} to get
\begin{align}\label{3.137}
\int_{0}^{d_1} |(\tilde{b}_{1},\tilde{b}_{2},\tilde{c})(x)|^2 dx\leq \frac{C}{(\s_0-\s)^2}\|e^{\s_0 x} \nu^{-1}wg\|^2_{L^\infty_{x,v}} d_1^2\, e^{-2\s d_1}.
\end{align}

Finally we consider the case for $\tilde{a}$. Multiplying \eqref{3.132} by $v_3\sqrt{\mu}$, we have that
\begin{align}\nonumber
\frac{d}{dx} \int_{\R^3} \tilde{f}(x,v)\cdot v_3^2 \sqrt{\mu} dv=0.
\end{align}
Integrating above equation over $[x,d]$ and using \eqref{3.119}, one obtain
\begin{align}\label{3.138}
\tilde{a}(x)=-2\tilde{c}(x)+\int_{\R^3}   (\mathbf{I-P})\tilde{f}(x,v)\cdot v_3^2 \sqrt{\mu} dv-\int_{\R^3}\tilde{f}(d_1,v)\cdot v_3\sqrt{\mu}  dv.
\end{align}
It follows from \eqref{3.138}, \eqref{3.134}, \eqref{3.137} and \eqref{3.129-0}, one has
\begin{align*}
\int_{0}^{d_1} |\tilde{a}(x)|^2 dx \leq \frac{C}{(\s_0-\s)}\|e^{\s_0 x}\nu^{-1} wg\|^2_{L^\infty_{x,v}} d_1^2e^{-2 \s d_1},
\end{align*}
which, together with \eqref{3.134} and \eqref{3.137}, yields that
\begin{align}\label{3.133}
\int_0^{d_1}\int_{\R^3}|\tilde{f}(x,v)|^2 dvdx\leq \frac{C}{(\s_0-\s)^2}\|e^{\s_0 x}\nu^{-1} wg\|^2_{L^\infty_{x,v}} d_1^2 e^{-2\s d_1}.
\end{align}

\vspace{1mm}

\noindent{\it Step 2. Convergence in $L^\infty$-norm.}
In the following, we shall use $t_{k}=t_{k}(t,x,v), X_{cl}(s;t,x,v), x_k=x_k(v,x)$ to be the back-time cycles defined for domain $[0,d_1]\times \R^3$. For later use, we denote $\tilde{h}:=w\tilde{f}$. Let $(x,v)\in [0,d_1]\times\R^3\backslash (\g_0\cup\g_{-})$, it follows from \eqref{3.132} that
\begin{align}\label{3.141}
\tilde{h}(x,v)&=e^{-\nu(v) (t-t_{k})} \tilde{h}(d_1,v_{k-1}) +\sum_{l=0}^{k-1} \int_{t_{l+1}}^{t_l}e^{-\nu(v) (t-s)}  K_w \tilde{h}(X_{cl}(s), v_l)ds,
\end{align}
with $k=1$ for $v_{0,3}<0$, and $k=2$ for $v_{0,3}>0$. We   will  use this summation convention in the following of this lemma.

From now on, we assume  $x\leq [0,\frac12 d_1]$. Then, if $v_{0,3}<0$, we have
\begin{align}\label{3.142}
t-t_k=t-t_1=\frac{(d_1-x) }{|v_{0,3}|}\geq \frac{d_1}{2|v_{0,3}|}.
\end{align}
If $v_{0,3}>0$, we obtain
\begin{align}\nonumber
t-t_k=t-t_2\geq t_1-t_2=\frac{d_1}{|v_{0,3}|}.
\end{align}
which, together with \eqref{3.142}, yields that
\begin{align}\label{3.143}
\nu(v) (t-t_k) \geq \f12\nu_0 d_1.
\end{align}
Hence, utilizing \eqref{3.143},  we always have
\begin{align}\label{3.144}
|e^{-\nu(v) (t-t_{k})} \tilde{h}(d_1,v_{k-1}) |
&\leq e^{-\f12 \nu_0 d_1} \Big(\|(w\bar{f}_{d_1},w\bar{f}_{d_2})\|_{L^\infty_{x,v}}+|w\bar{f}_{d_1}(d_1)|_{L^\infty(\g)}+|w\bar{f}_{d_2}(d_1)|_{L^\infty(\g)}\Big)\nonumber\\
&\leq \frac{C}{\s_0-\s}e^{-\nu_0 d_1}\|e^{\s_0x}\nu^{-1} wg\|_{L^\infty_{x,v}}.
\end{align}

For the second term on RHS of \eqref{3.141}, we use \eqref{3.141} again  to obtain
\begin{align}
&\left|\sum_{l=0}^{k-1} \int_{t_{l+1}}^{t_l}e^{-\nu(v) (t-s)}  K_w \tilde{h}(X_{cl}(s), v_l)ds\right|\nonumber\\
&=\sum_{l=0}^{k-1} \int_{t_{l+1}}^{t_l}e^{-\nu(v) (t-s)}  \int_{\R^3} |k_w(v_l, v') \tilde{h}(X_{cl}(s),v')|dv'ds\nonumber\\
&\leq \sum_{l=0}^{k-1} \int_{t_{l+1}}^{t_l}e^{-\nu(v) (t-s)}  \int_{\R^3} |k_w(v_l, v')| \nonumber\\ &\qquad\times\sum_{j=0}^{k-1}\int_{t'_{j+1}}^{t'_j}e^{-\nu(v') (s-s_1)}  \int_{\R^3} |k_w(v'_l, v'') \tilde{h}(X'_{cl}(s_1),v'')| dv'' ds_1 dv'ds\nonumber\\
&\quad+\frac{C}{\s_0-\s}e^{-\f12\nu_0 d_1}\|e^{\s_0x}\nu^{-1} wg\|_{L^\infty_{x,v}},\nonumber
\end{align}
where we have used \eqref{3.144} and denote $X'_{cl}(s_1)=X_{cl}(s_1;s,X_{cl}(s), v')$, $t'_j=t'_j(s_1;s, X_{cl}(s), v')$ and $v'_l$ to be the back-time cycle of $(s,X_{cl}(s), v')$. Then, by the same arguments as in Lemma \ref{lemS3.3}, we get
\begin{align}
&\|\tilde{h}\|_{L^\infty([0,\f12 d_1]\times \R^3)}+|\tilde{h}(0)|_{L^\infty(\gamma_+)}\nonumber\\
&\leq \f18 (\|\tilde{h}\|_{L^\infty([0,d_1]\times \R^3)}+|\tilde{h}(0)|_{L^\infty(\g_+)} ) + C(\|\bar{h}_{d_2}(d_1)\|_{L^\infty_v}+|\bar{h}_{d_1}(d_1)|_{L^\infty(\gamma)})\nonumber\\
&\quad+\frac{C}{\s_0-\s}e^{-\f12 \nu_0 d_1}\|e^{\s_0x} \nu^{-1}wg\|_{L^\infty_{x,v}}+C\|\tilde{f}\|_{L^2([0,d_1]\times\R^3)},\nonumber
\end{align}
which, together with \eqref{3.129-0} and \eqref{3.133}, yields that
\begin{align}\label{3.146}
&\|\tilde{h}\|_{L^\infty([0,\f12 d_1]\times \R^3)}+|\tilde{h}(0)|_{L^\infty(\gamma_+)}\nonumber\\
&\leq \frac{C}{\s_0-\s}\{e^{-\f12 \nu_0 d_1}+d_1^2 e^{-\frac12\s d_1}\}\|e^{\s_0x}\nu^{-1} wg\|_{L^\infty_{x,v}}\rightarrow0,\quad \mbox{as}\quad d_1\rightarrow\infty.
\end{align}
With the help of \eqref{3.146}, 
there exists a function $\mathfrak{f}(x,v)$ with $(x,v)\in\R_+\times \R^3$ so that $\|w(\bar{f}_{d}-\mathfrak{f})\|_{L^\infty([0,\f12 d]\times\R^3)}\rightarrow0$ as $d \rightarrow \infty$. The uniform bound \eqref{3.132-0}, \eqref{3.132-01}   follow from \eqref{3.129-0}, \eqref{3.123} and the strong convergence in $L^\infty_{x,v}$. It is  direct to see that $\mathfrak{f}(x,v)$ solves \eqref{1.7-2}. The continuity of $\mathfrak{f}$  follows directly from the $L^\infty_{x,v}$-convergence and the continuity of $\bar{f}_{d}$.

For the uniqueness, let $\mathfrak{f}_1,\mathfrak{f}_2$ be two solution of \eqref{1.7-1} with the bound \eqref{3.132-0} holds, then it holds that
\begin{align}\label{3.147}
\begin{cases}
v_3 \partial_x (\mathfrak{f}_1-\mathfrak{f}_2)+\FL (\mathfrak{f}_1-\mathfrak{f}_2)=0,\\
\mathfrak{f}_i(0,v)|_{v_3>0}=\mathfrak{f}_i(0,Rv),\ i=1,2,\\
\lim_{x\rightarrow\infty}\mathfrak{f}_i(x,v)=0,\ i=1,2.
\end{cases}
\end{align}
Multiplying \eqref{3.147} by $(\mathfrak{f}_1-\mathfrak{f}_2)$, it is direct to prove that
\begin{align}\nonumber
\int_0^\infty \|(\mathbf{I-P})(\mathfrak{f}_1-\mathfrak{f}_2)\|_{\nu}^2 dx=0.
\end{align}
That is, $(\mathfrak{f}_1-\mathfrak{f}_2)=\mathbf{P} (\mathfrak{f}_1-\mathfrak{f}_2)$. Then by the same arguments as in Lemma \ref{lem3.11} that
\begin{align}\nonumber
\int_0^\infty \|\mathbf{P}(\mathfrak{f}_1-\mathfrak{f}_2)\|_{L^2_v}^2 dx=0.
\end{align}
Thus, we prove $\mathfrak{f}_1\equiv\mathfrak{f}_2$. Therefore the proof of Lemma \ref{thm3.1} is completed. $\hfill\Box$

\

\section{Proof of Theorem \ref{thm1.1}} \label{M}
To prove the Theorem \ref{thm1.1}, we consider the following iterative sequence
\begin{align}\label{S3.123}
\begin{cases}
v_3\partial_x f^{j+1}+\FL f^{j+1}=\Gamma(f^j,f^j) +S,\\[1.5mm]
f^{j+1}(0,v)|_{v_3>0}= f^{j+1}(0,Rv)+f_b(Rv),\\[1.5mm]
\lim_{x\rightarrow\infty} f^{j+1}=0,
\end{cases}
\end{align}
for $j=0,1,2\cdots$ with $f^0\equiv0$.  It is direct to note  $\Gamma(f^j,f^j)\in \mathbb{N}^{\perp}$
and
\begin{align}\label{S3.126-1}
\|\nu^{-1} w \Gamma(f^j,f^j)\|_{L^\infty_v}\leq C \|wf^j\|^2_{L^\infty_v}.
\end{align}
Noting \eqref{S3.126-1}, and using Theorem  \ref{thm3.1}, we can solve \eqref{S3.123} inductively for $j=0,1,2,\cdots$. By taking $\frac34\s_0<\s<\s_0$, it follows from \eqref{S3.126-1} and \eqref{3.132-0} that
\begin{align}\label{S3.128}
&\|e^{\s x}wf^{j+1}\|_{L^\infty_{x,v}}+|wf^{j+1}(0)|_{L^\infty(\gamma)}\nonumber\\
&\leq \frac{C_1}{\s_0-\s}(|wf_b|_{L^\infty(\gamma)} +\|e^{\s_0 x}\nu^{-1}wS\|_{L^\infty_{x,v}})+\frac{C_1}{\s_0-\s} \|e^{\s x}wf^j\|^2_{L^\infty_{x,v}}.
\end{align}
We denote
\begin{equation*}
\delta=:|wf_b|_{L^\infty(\gamma)}+\|e^{\s_0 x} \nu^{-1} wS\|_{L^\infty_{x,v}}
\end{equation*}

By induction, we shall prove  that
\begin{equation}\label{S3.129}
\|e^{\s x} wf^{j}\|_{L^\infty_{x,v}}+|wf^{j}(0)|_{L^\infty(\gamma)}\leq \frac{2C_1\delta}{\s_0-\s},\quad\mbox{for} \  j=1,2,\cdots.
\end{equation}
Indeed, for $j=0$, it follows from $f^0\equiv0$ and \eqref{S3.128} that
\begin{equation*}
\|e^{\s x}wf^{1}\|_{L^\infty_{x,v}}+|wf^{1}(0)|_{L^\infty(\gamma)}\leq \frac{C_1\delta}{\s_0-\s}.
\end{equation*}
Now we assume that \eqref{S3.129} holds for $j=1,2\cdots, l$,  then we consider the case for $j=l+1$. Indeed it follows from \eqref{S3.128}  that
\begin{align*}
\|e^{\s x} wf^{l+1}\|_{L^\infty_{x,v}}+|wf^{l+1}(0)|_{L^\infty(\gamma)}&\leq \frac{C_1\delta}{\s_0-\s}+\frac{C_1}{\s_0-\s} \|e^{\s x}wf^l\|^2_{L^\infty_{x,v}}\nonumber\\
&\leq \frac{C_1\delta}{\s_0-\s}\Big(1+4(\frac{C_1}{\s_0-\s})^2\delta\Big)\leq  \frac{\f{3}{2}C_1\delta}{\s_0-\s},
\end{align*}
where we have used  $\eqref{S3.129}$ with $j=l$, and chosen $\delta\leq \d_0$ with $\delta_0$ small enough such that  $4 (\frac{C_1}{\s_0-\s})^2 \delta_0\leq 1/2$.
Therefore we have proved \eqref{S3.129} by induction.

Finally we consider the convergence of sequence  $f^j$. For the difference $f^{j+1}-f^j$, we have
\begin{equation}\label{S3.131}
\begin{cases}
\dis v_3\partial_x (f^{j+1}-f^j)+\FL (f^{j+1}-f^j)=\Gamma(f^{j}-f^{j-1},f^{j})+\Gamma(f^{j-1},f^{j}-f^{j-1}),\\[2mm]
\dis (f^{j+1}-f^j)(0,v)|_{v_3>0}= (f^{j+1}-f^j)(0,Rv),\\[2mm]
\lim_{x\rightarrow\infty}(f^{j+1}-f^j)(x,v)=0.
\end{cases}
\end{equation}
Applying \eqref{3.132-0} to \eqref{S3.131}, we have
\begin{align}\label{S3.133}
&\|e^{\s x}w\{f^{j+1}-f^j\}\|_{L^\infty_{x,v}}+|w\{f^{j+1}-f^j\}(0)|_{L^\infty}\nonumber\\
&\leq \frac{C}{\s_0-\s}\Big\{ \|e^{\s_0 x}\nu^{-1}w\Gamma(f^{j}-f^{j-1},f^{j})\|_{L^\infty_{x,v}} +\|e^{\s_0 x}\nu^{-1}w\Gamma(f^{j-1},f^{j}-f^{j-1})\|_{L^\infty_{x,v}}\Big\}\nonumber\\
&\leq \frac{C}{\s_0-\s}[\|e^{\s x}wf^{j}\|_{L^\infty_{x,v}}+\|e^{\s x}wf^{j-1}\|_{L^\infty_{x,v}}]\cdot \|e^{\s x}w(f^{j}-f^{j-1})\|_{L^\infty_{x,v}}\nonumber\\
&\leq \frac{C\delta_0}{(\s_0-\s)^2} \|e^{\s x}w(f^{j}-f^{j-1})\|_{L^\infty_{x,v}} 
\leq \frac12 \|e^{\s x}w(f^{j}-f^{j-1})\|_{L^\infty_{x,v}},
\end{align}
where we have used \eqref{S3.129} and further taken $\delta_0>0$ small such that $\dis \frac{C}{(\s_0-\s)^2}\delta_0\leq 1/2$. Hence $f^j$ is a Cauchy sequence in $L^\infty_{x,v}$, then we obtain the solution by taking the limit $f=\lim_{j\rightarrow\infty} f^j$.  The uniqueness can also be obtained by using the inequality as \eqref{S3.133}. The continuity of $f$ is   a direct consequence of $L^\infty_{x,v}$-convergence. 
Therefore we complete the proof of Theorem \ref{thm1.1}. \qed



\

\noindent{\bf Acknowledgments.} 
Feimin Huang's research is partially supported by National Natural Sciences Foundation of China No. 11688101.   Zaihong Jiang partially supported by  ZJNSF Grant No.LY19A010016. Yong Wang's research  is partially supported by National Natural Sciences Foundation of China No. 11771429, 11671237 and  11688101.

\end{document}